\documentclass[11pt,a4paper]{article}
\usepackage[leqno]{amsmath}
\usepackage{amssymb,amsthm}
\usepackage[a4paper,textwidth=150mm,textheight=245mm,centering]{geometry}

\newtheorem{thm}{Theorem}[section]
\newtheorem{lem}[thm]{Lemma}
\numberwithin{equation}{section}

\allowdisplaybreaks

\begin{document}

\title{The sum of sizes of nontrivial cross-intersecting families}
\date{}
\author{Huajun Zhang\\
Shaoxing University\\ Shaoxing 312000, P. R. China\\
\texttt{huajunzhang@usx.edu.cn}}
\maketitle
\begin{abstract}
Families $\mathcal A_1,\ldots,\mathcal A_m\subseteq\binom{[n]}k$
are cross-intersecting if every member of one family meets every
member of each other family. A family is nontrivial if its members
have no common element.

For $n>2k$ and $m\ge3$, we determine the maximum sum of the sizes
of $m$ nonempty, nontrivial cross-intersecting families. We prove
the sharp bound
\[
\sum_{i=1}^{m}|\mathcal A_i|
\le \max\bigl\{m\,h(n,k,1),\ g(n,k)+3(m-1)\bigr\},
\]
where
\[
h(n,k,1)=\binom{n-1}{k-1}-\binom{n-k-1}{k-1}+1
\]
and
\[
g(n,k)=\binom nk-2\binom{n-k}k
       -\binom{n-k-1}{k-1}+\binom{n-2k+1}k.
\]
We also characterize all equality cases. Two types of extremal
configurations arise: copies of a maximum nontrivial intersecting
family, and copies of a specific three-member nontrivial
intersecting family together with its full cross-intersecting
companion. The proof applies a known restricted lexicographic
comparison to obtain the weighted bound for shifted families,
and uses local generator reductions for critically shifted
families. These reductions exploit shiftedness
after removing the common core of a section, while preserving
nontriviality and the required shadow relations between sections.
\end{abstract}

\section{Introduction}

Let $[n]=\{1,2,\ldots,n\}$, and write $2^{[n]}$ for its power set
and $\binom{[n]}k$ for the family of all $k$-element subsets of $[n]$.
Families
$\mathcal F_1,\ldots,\mathcal F_m\subseteq\binom{[n]}k$ are
\emph{cross-intersecting} if $F\cap F'\ne\emptyset$ whenever
$F\in\mathcal F_i$, $F'\in\mathcal F_j$, and $i\ne j$.
We study the sum $|\mathcal F_1|+\cdots+|\mathcal F_m|$ when every family
is nonempty and nontrivial, meaning that its members have no common element.

Cross-intersecting families generalize intersecting families: a family
$\mathcal F$ is \emph{intersecting} if every two of its members meet.
The classical theorem of Erd\H{o}s, Ko and Rado \cite{EKR} determines
the maximum size of an intersecting family.

\begin{thm}[Erd\H{o}s--Ko--Rado \cite{EKR}]
Let $n\ge2k$, and let $\mathcal F\subseteq\binom{[n]}k$ be intersecting.
Then
\[
|\mathcal F|\le\binom{n-1}{k-1}.
\]
Moreover, if $n>2k$, equality holds if and only if
\[
\mathcal F=\left\{F\in\binom{[n]}k:x\in F\right\}
\]
for some $x\in[n]$.
\end{thm}

The extremal family in this theorem is called a \emph{star}.
For $n>2k$, requiring $\mathcal F$ to be nontrivial excludes every
maximum intersecting family. Hilton and Milner \cite{hm} determined
the maximum under this additional requirement.

\begin{thm}[Hilton--Milner \cite{hm}]
Let $n>2k$, and let $\mathcal F\subseteq\binom{[n]}k$ be a nontrivial
intersecting family. Then
\begin{equation}\label{hm1}
|\mathcal F|\le\binom{n-1}{k-1}-\binom{n-k-1}{k-1}+1.
\end{equation}
\end{thm}

Write $[i,j]=\{i,i+1,\ldots,j\}$ for $i\le j$, with $[0]=\emptyset$.
For $k\ge2$, define
\[
\mathcal H(n,k,1)=
\left\{A\in\binom{[n]}k:1\in A,\ A\cap[2,k+1]\ne\emptyset\right\}
\cup\{[2,k+1]\},
\]
and put
\[
h(n,k,1)=|\mathcal H(n,k,1)|
=\binom{n-1}{k-1}-\binom{n-k-1}{k-1}+1.
\]
Equality holds in \eqref{hm1} if and only if, up to isomorphism,
$\mathcal F=\mathcal H(n,k,1)$ or, when $k=3$,
\[
\mathcal F=\left\{F\in\binom{[n]}3:|F\cap[3]|\ge2\right\}.
\]

For cross-intersecting families, Hilton \cite{hilton} determined
the maximum sum when empty families are allowed.

\begin{thm}
Let $\mathcal A_1,\ldots,\mathcal A_m\subseteq\binom{[n]}k$ be
cross-intersecting families with $\mathcal A_1\ne\emptyset$.
If $n\ge2k$, then
\[
\sum_{i=1}^m|\mathcal A_i|
\le
\begin{cases}
\binom nk, & \text{if }m\le n/k,\\[4pt]
m\binom{n-1}{k-1}, & \text{if }m\ge n/k.
\end{cases}
\]
\end{thm}

The two values are attained by one full family and $m-1$ empty
families, and by $m$ copies of a star, respectively.
Borg \cite{pborg1} gave a short proof and studied related extensions
to labeled sets \cite{borg} and permutations \cite{pborg}.
Wang and the author \cite{wzhj} extended the result to symmetric
systems, including finite sets, finite vector spaces, and permutations.

Requiring every family to be nonempty changes the problem.
For two families of possibly different uniformities, the results
of Hilton and Milner \cite{hm} and Frankl and Tokushige \cite{FT}
give
\[
|\mathcal A|+|\mathcal B|
\le\binom n\ell-\binom{n-k}\ell+1
\]
whenever $\mathcal A\subseteq\binom{[n]}k$ and
$\mathcal B\subseteq\binom{[n]}\ell$ are nonempty and
cross-intersecting, $n\ge k+\ell$, and $k\le\ell$.
This bound is attained by a single $k$-set and all $\ell$-sets
meeting it. Wang and the author \cite{wz2} extended this inequality
to a setting that also includes finite vector spaces and permutations.

Nontriviality excludes this singleton construction as well as stars.
For two nontrivial cross-intersecting families, M\"ors \cite{mh}
obtained the following threshold bound.
The formulation below, including its strictness assertion, also
follows from F\"uredi's theorem \cite[Theorem, p.~333]{furedi95}.

\begin{thm}[M\"ors \cite{mh}]\label{mh}
Let $n,k,\ell$ be positive integers with $n>k+\ell$.
Let $\mathcal A\subseteq\binom{[n]}k$ and
$\mathcal B\subseteq\binom{[n]}\ell$ be nonempty, nontrivial
cross-intersecting families. Then at least one of the following
inequalities holds:
\[
|\mathcal A|\le
\binom{n-1}{k-1}-\binom{n-\ell-1}{k-1}+1,
\]
\[
|\mathcal B|\le
\binom{n-1}{\ell-1}-\binom{n-k-1}{\ell-1}+1.
\]
Moreover, if one of these inequalities fails, the other holds strictly.
\end{thm}

For equal uniformities, this implies that at most one family in a
nontrivial cross-intersecting tuple can have size greater than
$h(n,k,1)$. It does not, however, determine the maximum sum.
Frankl \cite{fr24} settled the sum problem for two families.

\begin{thm}[\cite{fr24}]\label{cab1}
Let $n>2k$, and let $\mathcal A,\mathcal B\subseteq\binom{[n]}k$
be nonempty, nontrivial cross-intersecting families. Then
\[
|\mathcal A|+|\mathcal B|
\le\binom nk-2\binom{n-k}k+\binom{n-2k}k+2.
\]
Moreover, if $k\ge3$, equality holds if and only if, after
interchanging $\mathcal A$ and $\mathcal B$ if necessary,
\[
\mathcal A=\{A_1,A_2\},\qquad
\mathcal B=
\left\{B\in\binom{[n]}k:
B\cap A_1\ne\emptyset,\ B\cap A_2\ne\emptyset\right\},
\]
for some disjoint $A_1,A_2\in\binom{[n]}k$.
\end{thm}

For three or more families, the two-family extremal construction
does not extend by simply repeating its smaller family:
a family consisting of two disjoint sets is not cross-intersecting
with a second copy of itself. Repeated families must themselves be
intersecting. This brings the Hilton--Milner bound back into the
sum problem and suggests a comparison between large nontrivial
intersecting families and small ones with large cross-intersecting
companions.

\subsection{Main result}

One construction consists of $m$ copies of a maximum nontrivial
intersecting family and has total size $m\,h(n,k,1)$.
For the second construction, with $k\ge2$, define
\[
\mathcal M=
\left\{[k],\ [k,2k-1],\ [k+1]\setminus\{k\}\right\},
\]
and let
\[
\mathcal N=
\left\{A\in\binom{[n]}k:
A\cap B\ne\emptyset\text{ for every }B\in\mathcal M\right\}.
\]
The family $\mathcal M$ is nontrivial and intersecting.
Since $\mathcal M\subseteq\mathcal N$, the family $\mathcal N$
is also nontrivial. Thus
$\mathcal N,\mathcal M,\ldots,\mathcal M$ is an admissible
cross-intersecting tuple. Put
\[
g(n,k)=|\mathcal N|
=\binom nk-2\binom{n-k}k-\binom{n-k-1}{k-1}
+\binom{n-2k+1}k.
\]
The total size of this tuple is $g(n,k)+3(m-1)$.

Our main result shows that these two constructions determine
the maximum sum and account for all equality cases.

\begin{thm}\label{mth}
Let $m,n,k$ be positive integers with $m\ge3$ and $n>2k$.
Let $\mathcal A_1,\ldots,\mathcal A_m\subseteq\binom{[n]}k$
be nonempty, nontrivial cross-intersecting families. Then
\[
\sum_{i=1}^m|\mathcal A_i|
\le\max\left\{m\,h(n,k,1),\ g(n,k)+3(m-1)\right\}.
\]
Moreover, equality holds if and only if, up to a permutation of
$[n]$ and a relabeling of the families, one of the following occurs:
\begin{itemize}
\item[\textup{(i)}]
$\mathcal A_1=\cdots=\mathcal A_m=\mathcal F$, where $\mathcal F$
is a maximum nontrivial intersecting family in $\binom{[n]}k$, and
\[
m\,h(n,k,1)\ge g(n,k)+3(m-1).
\]
\item[\textup{(ii)}]
$\mathcal A_1=\mathcal N$ and
$\mathcal A_2=\cdots=\mathcal A_m=\mathcal M$, and
\[
m\,h(n,k,1)\le g(n,k)+3(m-1).
\]
\end{itemize}
\end{thm}

\subsection{Main ideas}

We first order the families by size and compare their sum with
$(m-1)|\mathcal A_2|+|\mathcal A_1|$. This leads to a weighted
problem for a cross-intersecting pair. The main difficulty is
to simplify the pair while preserving nontriviality.

Shifting preserves the sizes of families and cross-intersection,
but it may create a common element. We therefore apply simultaneous
shifts only while the required nontriviality is preserved.
The relevant family eventually becomes either shifted or
\emph{critically shifted}: it is not shifted, but every shift
that changes it makes it trivial. In the shifted case, the
restricted lexicographic comparison of Wu and Zhang
\cite{wuzhang26} reduces the problem to initial segments.
We then prove the required weighted estimate, using their
threshold classification to recover equality.

The critically shifted case requires local modifications that
respect the relations between sections. We choose a critical pair
of coordinates and consider the sections containing exactly one
or both of these coordinates. After removing the common core
of a single-index section, the remaining family is shifted on
its residual ground set. This local shiftedness allows us to
modify generators without requiring the whole family to be shifted.

We contract the generators of a single-index section while keeping
the double section fixed. The contraction preserves the required
relative shadow inclusion and the nontriviality of the reconstructed
family. It removes some members from one family and allows new
members in its cross-intersecting companion. By counting these
new members directly, we obtain strict improvements in the weighted
sum that exclude intermediate section structures.
The remaining cases are handled by explicit estimates, and the
equality cases are recovered by tracing the reductions and shifts.

Section~2 collects the auxiliary inequalities.
Section~3 introduces generators and records the neighborhood
comparison results used below. Section~4 proves the weighted bound for
shifted families. Section~5 develops the reductions for critically
shifted families and proves Theorem~\ref{mth}.

\section{Auxiliary inequalities}

This section collects the comparisons used to rule out the other
configurations that arise in the proof. We first compare weighted
binomial expressions, then estimate the neighborhood of two
cross-intersecting disjoint pairs. The final lemma gives two
additional comparisons for the critically shifted case.

We use $\binom ab=0$ whenever $b<0$ or $a<b$, including negative
upper arguments; otherwise $\binom ab$ has its usual combinatorial
meaning. Empty sums are $0$ and empty products are $1$.
For a family $\mathcal F$ of subsets of $[n]$, define
\[
N_\ell(\mathcal F)=
\left\{B\in\binom{[n]}\ell:
B\cap F=\emptyset\text{ for some }F\in\mathcal F\right\}.
\]
Its complement
\[
\overline N_\ell(\mathcal F)=\binom{[n]}\ell\setminus N_\ell(\mathcal F)
\]
is the full $\ell$-uniform family cross-intersecting with $\mathcal F$.
We write $N_\ell(F)$ for $N_\ell(\{F\})$.
When the ground set is smaller than $[n]$, we specify it explicitly.

Our first lemma transfers a weighted binomial comparison to two
related comparisons. The condition $q>i+j$ makes both conclusions strict.

\begin{lem}\label{mpij}
Let $i,j,q$ be positive integers with $q>i+j$, and let $x,y$ be
positive real numbers. If
\[
x\binom qi\le y\binom qj,
\]
then
\begin{equation}\label{eqij1}
x\binom{q-1}{i-1}<y\binom{q-1}{j},
\end{equation}
and
\begin{equation}\label{eqij2}
x\binom q{i-1}<y\binom q{j+1}.
\end{equation}
\end{lem}

\begin{proof}
The first comparison follows from
\[
\frac{x\binom{q-1}{i-1}}{y\binom{q-1}{j}}
=\frac{i}{q-j}\frac{x\binom qi}{y\binom qj}
\le\frac{i}{q-j}<1.
\]
For the second, we similarly have
\[
\frac{x\binom q{i-1}}{y\binom q{j+1}}
=\frac{i(j+1)}{(q-i+1)(q-j)}
 \frac{x\binom qi}{y\binom qj}
\le\frac{i(j+1)}{(q-i+1)(q-j)}<1.
\]
The last inequality holds because
$(q-i+1)(q-j)-i(j+1)=(q+1)(q-i-j)>0$.
\end{proof}

The next two lemmas compare the value obtained from
$\binom{[k+1]}k$ and its full cross-intersecting family with the
candidate extremal values. They will exclude this clique
construction in the shifted case. We begin with different uniformities.

\begin{lem}\label{inma}
Let $m,n,k,\ell$ be positive integers with $k\ge2$, $n>k+\ell$,
and $m\binom nk<\binom n\ell$. Then
\[
\binom n\ell-2\binom{n-k}{\ell}+\binom{n-2k}{\ell}+2m
>\sum_{i=2}^{k+1}\binom{k+1}{i}\binom{n-k-1}{\ell-i}+m(k+1).
\]
\end{lem}

\begin{proof}
The assumptions force $\ell>k$: if $\ell\le k$ and $n>k+\ell$,
then $\binom n\ell\le\binom nk$, contrary to the hypothesis.
We first bound $m$. Since $n\ge k+\ell+1$,
\[
\frac{\binom{n-k-2}{\ell-2}\binom nk}{\binom n\ell}
=\frac{\ell(\ell-1)}{k!}
 \frac{n-\ell}{n-k}\frac{n-\ell-1}{n-k-1}
 \prod_{j=2}^{k-1}(n-\ell-j)
\ge\frac{(k+1)(\ell-1)}{\ell+1}>1.
\]
Each factor depending on $n$ is nondecreasing, so the first
inequality follows by setting $n=k+\ell+1$; the product is $1$
when $k=2$. Thus
\begin{equation}\label{mnk}
m<\frac{\binom n\ell}{\binom nk}<\binom{n-k-2}{\ell-2}.
\end{equation}

Let $\Delta$ be the left-hand side minus the right-hand side
of the asserted inequality. Then
\[
\Delta=(k-1)\binom{n-k-1}{\ell-1}
-\sum_{i=2}^{k}\binom{n-k-i}{\ell-1}-(k-1)m.
\]
Each summand is at most $\binom{n-k-2}{\ell-1}$, so
\[
\Delta\ge(k-1)\left(\binom{n-k-2}{\ell-2}-m\right)>0,
\]
where the strict inequality is \eqref{mnk}.
\end{proof}

For equal uniformities, we compare with the Hilton--Milner value
when the weight is large and with the three-member construction
when it is small.

\begin{lem}\label{mklc}
Let $n,k,m$ be positive integers with $k\ge3$, $n>2k$, and $m\ge2$.
Then
\[
\sum_{i=2}^{k}\binom{k+1}{i}\binom{n-k-1}{k-i}+m(k+1)
<\max\left\{(m+1)h(n,k,1),\,g(n,k)+3m\right\}.
\]
\end{lem}

\begin{proof}
Write $S$ for the sum on the left before the term $m(k+1)$.
We have
\[
S=\sum_{i=1}^{k}\binom{n-i}{k-1}
-k\binom{n-k-1}{k-1}.
\]
Also,
\[
h(n,k,1)=1+\sum_{i=2}^{k+1}\binom{n-i}{k-2}>k+1.
\]

Suppose first that $m\ge k-1$.
Since $h(n,k,1)>k+1$, the difference
$(m+1)h(n,k,1)-S-m(k+1)$ is increasing in $m$.
It is therefore enough to take $m=k-1$, when
\begin{align*}
kh(n,k,1)-S-(k-1)(k+1)
&=(k-1)\binom{n-1}{k-1}
  -\sum_{i=2}^{k}\binom{n-i}{k-1}-k(k-1)+1\\
&\ge(k-1)\left(\binom{n-2}{k-2}-k\right)+1>0.
\end{align*}
Here we used $\binom{n-2}{k-2}\ge\binom{2k-1}{k-2}>k$.

Now suppose that $2\le m\le k-2$, so $k\ge4$.
The definition of $g(n,k)$ gives
\begin{align*}
g(n,k)+3m-S-m(k+1)
&=(k-2)\binom{n-k-1}{k-1}
  -\sum_{i=2}^{k-1}\binom{n-k-i}{k-1}-m(k-2)\\
&\ge(k-2)\left(\binom{n-k-2}{k-2}-m\right)>0,
\end{align*}
since $\binom{n-k-2}{k-2}\ge\binom{k-1}{k-2}=k-1>m$.
This proves the lemma.
\end{proof}

We next consider two cross-intersecting pairs of disjoint sets.
The following neighborhood estimate will exclude configurations
in which two families each consist of such a pair.

\begin{lem}\label{abc}
Let $A_1,A_2,B_1,B_2\in\binom{[n]}k$ satisfy
$A_1\cap A_2=B_1\cap B_2=\emptyset$ and
$A_i\cap B_j\ne\emptyset$ for $i,j\in[2]$.
If $n>k+\ell$ and $\ell\ge2$, then
\[
|N_\ell(\{A_1,A_2,B_1,B_2\})|
>2\binom{n-k}{\ell}+\binom{n-k-1}{\ell-1}
-\binom{n-2k+1}{\ell}.
\]
\end{lem}

\begin{proof}
The assumptions imply $k\ge2$.
Since $A_1\cap A_2=\emptyset$,
\[
|N_\ell(\{A_1,A_2\})|
=2\binom{n-k}{\ell}-\binom{n-2k}{\ell}.
\]
It remains to show that
\begin{equation}\label{abc-extra-neighborhood}
|N_\ell(\{B_1,B_2\})\setminus N_\ell(\{A_1,A_2\})|
>\binom{n-k-1}{\ell-1}-\binom{n-2k}{\ell-1}.
\end{equation}

For $i\in[2]$, put
\[
C_i=A_1\cap B_i,\qquad D_i=A_2\cap B_i.
\]
These sets are nonempty and disjoint, with $|C_i|+|D_i|\le k$.
Let $E_i=N_\ell(B_i)\setminus N_\ell(\{A_1,A_2\})$.
Thus $E_i$ consists of the $\ell$-sets that avoid $B_i$ and meet
both $A_1$ and $A_2$. Hence
\begin{align*}
|E_i|
&=\binom{n-k}{\ell}-\binom{n-2k+|C_i|}{\ell}
  -\binom{n-2k+|D_i|}{\ell}
  +\binom{n-3k+|C_i|+|D_i|}{\ell}\\
&\ge\binom{n-k}{\ell}-\binom{n-2k+|C_i|}{\ell}
  -\binom{n-k-|C_i|}{\ell}+\binom{n-2k}{\ell}\\
&\ge\binom{n-k-1}{\ell-1}-\binom{n-2k}{\ell-1}.
\end{align*}
For the first inequality, work on $[n]\setminus B_i$.
A set must meet both $A_1\setminus B_i$ and $A_2\setminus B_i$,
which are disjoint. Restricting the latter to $|C_i|$ elements
can only decrease the count. For the second inequality, discrete
convexity gives, for $1\le|C_i|\le k-1$,
\[
\binom{n-2k+|C_i|}{\ell}+\binom{n-k-|C_i|}{\ell}
\le\binom{n-2k+1}{\ell}+\binom{n-k-1}{\ell}.
\]

We also need the equality conditions in this estimate.
If $2\le|C_i|\le k-2$, then $k\ge4$, and the convexity
inequality is strict: the drop from an endpoint to the next
value is
\[
\binom{n-k-2}{\ell-1}-\binom{n-2k+1}{\ell-1}>0,
\]
because $n-k-2\ge\ell-1$ and $k\ge4$.
If $|C_i|=1$ and $|D_i|<k-1$, the first inequality is strict.
Indeed, after restricting $A_2\setminus B_i$ to a singleton,
an $\ell$-set can avoid that singleton and meet both
$A_1\setminus B_i$ and the removed part; it can be completed
because $n-k-1\ge\ell$.
Finally, $|C_i|=k-1$ forces $|D_i|=1$.
Thus equality holds in the bound for $|E_i|$ exactly when
$\{|C_i|,|D_i|\}=\{1,k-1\}$.

If either of these two bounds is strict,
\eqref{abc-extra-neighborhood} follows at once.
Otherwise, both $B_1$ and $B_2$ lie in $A_1\cup A_2$.
Their disjointness then gives $B_1\cup B_2=A_1\cup A_2$.
After interchanging $B_1,B_2$ and relabeling, we may assume that
\[
A_1=[k],\quad A_2=[k+1,2k],\quad
B_1=[k+1]\setminus\{k\},\quad B_2=[k,2k]\setminus\{k+1\}.
\]

Every member of $E_1$ contains $k$ and avoids $k+1$,
whereas every member of $E_2$ contains $k+1$ and avoids $k$.
Hence $E_1$ and $E_2$ are disjoint, and
\[
|E_1\cup E_2|
=2\left(\binom{n-k-1}{\ell-1}-\binom{n-2k}{\ell-1}\right)
>\binom{n-k-1}{\ell-1}-\binom{n-2k}{\ell-1}.
\]
The strict inequality follows from $n-k-1\ge\ell$,
$k\ge2$, and $\ell\ge2$.
Since $E_1\cup E_2$ is the neighborhood difference in
\eqref{abc-extra-neighborhood}, the proof is complete.
\end{proof}

The last lemma supplies two comparisons for the critical-shifting
arguments. Part~\textup{(i)} records the bound on $m$ from
\eqref{mnk}. Part~\textup{(ii)} compares a second difference of
binomial coefficients with a weighted binomial coefficient.

\begin{lem}\label{m2n-binomial}
Let $n,m,k,\ell$ be positive integers satisfying $m\ge2$,
$n>k+\ell>2k$, $k\ge2$, and $m\binom nk<\binom n\ell$.
Then the following hold.
\begin{enumerate}
\item[(i)] $m<\binom{n-k-2}{\ell-2}$.
\item[(ii)] If $k\ge3$ and $2\le t\le k-1$, then
\[
\binom{n-t}{\ell}-2\binom{n-k}{\ell}
+\binom{n-2k+t}{\ell}+2m
>m\binom{n-t}{k-t}.
\]
\end{enumerate}
\end{lem}

\begin{proof}
Part~\textup{(i)} is \eqref{mnk}.
For part~\textup{(ii)}, put $R(n)=\binom n\ell/\binom nk$,
so $m<R(n)$, and let $D_t(n)$ denote the binomial expression
on the left before the term $2m$.
We normalize by $R(n)$ and compare at the smallest admissible
value of $n$, namely $k+\ell+1$.
The assumptions also give $\ell\ge k+1$.

First suppose that $t=k-1$.
Here $D_{k-1}(n)=\binom{n-k-1}{\ell-2}$.
The ratio $D_{k-1}(n)/((n-k-1)R(n))$ is strictly increasing
in $n$: the ratio of its values at $n+1$ and $n$ is
\[
1+\frac{(n-k)(k-3)+2\ell-k-3}
        {(n-k-\ell+2)(n-k+1)}>1.
\]
At $n=k+\ell+1$, its value is
$\frac{(k+1)(\ell-1)}{2(\ell+1)}\ge\frac65>1$.
Consequently,
\[
D_{k-1}(n)+2m>m(n-k-1)+2m
=m\binom{n-k+1}{1},
\]
as required.

Now suppose that $2\le t\le k-2$.
We have
\[
\frac{D_t(n)}{R(n)\binom{n-t}{k-t}}
=\sum_{u=1}^{k-t}\sum_{v=1}^{k-t}
 \frac{\binom{n-t-u-v}{\ell-2}}
      {R(n)\binom{n-t}{k-t}}.
\]
Each summand is nondecreasing in $n$.
When positive, the ratio of its values at $n+1$ and $n$ is
\[
1+\frac{(t-2)(n-t+1)+(u+v)(\ell-t)}
        {(n-t-u-v-\ell+3)(n-t+1)}>1;
\]
a zero summand remains zero or becomes positive.
At $n=k+\ell+1$, the normalized sum is
\[
\frac{k+1}{k-t+1}
\left(1-\frac{2(\ell+1)}
{\binom{k+\ell-t+1}{k-t+1}}\right)>1,
\]
where the inequality follows from
\[
t\binom{k+\ell-t+1}{k-t+1}
\ge2\binom{\ell+3}{3}>2\ell(\ell+1)
\ge2(k+1)(\ell+1).
\]
It follows that
\[
D_t(n)>R(n)\binom{n-t}{k-t}>m\binom{n-t}{k-t},
\]
which proves part~\textup{(ii)}.
\end{proof}

\section{Generators, lexicographic order and neighborhoods}

We introduce the generator notation used in the local reductions
and recall the lexicographic neighborhood estimates needed below.
The restricted comparison and its threshold equality classification
are taken from Wu and Zhang \cite{wuzhang26}.

We use the standard shifting notation; see \cite{F87}.
For $i<j$, the shift $s_{ij}(\mathcal F)$ replaces $F\in\mathcal F$
by $(F\setminus\{j\})\cup\{i\}$ whenever $j\in F$, $i\notin F$,
and the latter set does not already belong to $\mathcal F$.
A family is \emph{shifted} if every such shift fixes it.
A nontrivial family is \emph{critically shifted} if it is not
shifted and every shift that changes it makes it trivial.

\subsection{Generators}

For $0\leq r\leq n$ and
$\Omega\subseteq\bigcup_{i=0}^{r}\binom{[n]}i$, define
\[
\langle\Omega\rangle_r
=\left\{F\in\binom{[n]}r:
G\subseteq F\text{ for some }G\in\Omega\right\}.
\]
We call $\Omega$ a \emph{generating family} of $\mathcal F$ if
$\mathcal F=\langle\Omega\rangle_r$, and call its members
\emph{generators}. A generator containing another generator is
redundant and may be deleted. We use the following canonical choice:
\[
\mathcal G_r(\mathcal F)=\min_{\subseteq}
\left\{G\subseteq[n]:|G|\leq r,\;
\{F\in\tbinom{[n]}r:G\subseteq F\}\subseteq\mathcal F\right\},
\]
where $\min_{\subseteq}$ denotes the inclusion-minimal members.
Then $\mathcal F=\langle\mathcal G_r(\mathcal F)\rangle_r$.
The empty family has no generators, and the full level is generated
by $\emptyset$. We use the same notation on smaller ground sets,
which will be specified when needed.

The next lemma expresses cross-intersection and common intersections
in terms of generators.

\begin{lem}\label{generator-rules}
Let $n\geq k+\ell$, $\mathcal A=\langle\Omega\rangle_k$ and
$\mathcal B=\langle\Lambda\rangle_\ell$.
\begin{enumerate}
\item[\textup{(i)}] The families $\mathcal A,\mathcal B$ are cross-intersecting if and only if
$G\cap H\neq\emptyset$ for every $G\in\Omega$ and $H\in\Lambda$.
Equivalently,
\[
\overline N_\ell(\mathcal A)
=\{B\in\tbinom{[n]}\ell:
B\cap G\neq\emptyset\text{ for all }G\in\Omega\}.
\]
\item[\textup{(ii)}] If $n>k$ and $\mathcal A\neq\emptyset$, then
$\bigcap\mathcal A=\bigcap\Omega$.
\end{enumerate}
\end{lem}

\begin{proof}
For \textup{(i)}, intersecting generators give intersecting
supersets. Conversely, if $G\cap H=\emptyset$, extend $G$ to a
$k$-set avoiding $H$, then extend $H$ to an $\ell$-set avoiding
that $k$-set. Both extensions are possible because $n\geq k+\ell$.
The description of $\overline N_\ell(\mathcal A)$ follows in the
same way, by extending a generator disjoint from a prescribed
$\ell$-set.

For \textup{(ii)}, every member of $\mathcal A$ contains a generator,
so $\bigcap\Omega\subseteq\bigcap\mathcal A$. If a point is absent
from some generator, then $n>k$ allows a $k$-set containing that
generator and avoiding the point. This proves the reverse inclusion.
\end{proof}

\subsection{Lexicographic order}

For distinct $F,G\in\binom{[n]}k$, define the lexicographic order by
\[
F\prec G\quad\Longleftrightarrow\quad
\min(F\setminus G)<\min(G\setminus F).
\]
Let $\mathcal L([n],k,m)$ denote the first $m$ members in this order.
For $n>k$, every integer $1\le m\le\binom nk$ has a unique canonical
expansion
\[
m=\sum_{i=0}^{s}\binom{n-a_i}{k-a_i+i},
\]
where $0\le a_0<\cdots<a_s$, $a_i\le k+i$, and $0\le s<n-k$.
The co-Kruskal--Katona theorem states that an initial lexicographic
segment minimizes the disjointness neighborhood among families
of a given size; see \cite[Theorem~2.3]{wuzhang26}.
In terms of the canonical expansion, it takes the following form.

\begin{thm}[co-Kruskal--Katona theorem]\label{cokk}
Let $n, k, \ell \in \mathbb{Z}^+$ with $n \geq k + \ell$. If $\mathcal{A} \subseteq \binom{[n]}{k}$ satisfies
\[
|\mathcal{A}| = \sum_{i=0}^s \binom{n - a_i}{k - a_i + i}
\]
in its canonical expansion, then
\[
|N_\ell(\mathcal{A})| \geq \sum_{i=0}^s \binom{n - a_i}{\ell-i}.
\]
Equality holds if $\mathcal{A} =\mathcal{L}([n],k,|\mathcal{A}|)$.
\end{thm}

The corresponding replacement principle for cross-intersecting
pairs is Hilton's lemma.

\begin{lem}[Hilton \cite{Hil1976}] \label{lem1}
Let $k, \ell, n$ be positive integers. If $\mathcal{A} \subseteq \binom{[n]}{k}$ and $\mathcal{B} \subseteq \binom{[n]}{\ell}$ are
cross-intersecting, then the initial segments $\mathcal{L}([n], k, |\mathcal{A}|)$ and $\mathcal{L}([n], \ell, |\mathcal{B}|)$ are also
cross-intersecting.
\end{lem}

\subsection{Restricted neighborhood estimates}

Lexicographic replacement preserves cross-intersection but may
destroy nontriviality. For $n>k+\ell$ and $k,\ell\ge2$, put
\[
\mathcal B(n,k,\ell)
=\left\{A\in\binom{[n]}k:|A\cap[\ell+1]|\ge2\right\},
\]
and let $\mathcal L_\ell([n],k,m)$ be its first $m$ members in
lexicographic order, where $0\le m\le|\mathcal B(n,k,\ell)|$.
By \cite[Lemma~2.5]{wuzhang26}, a shifted family with nonempty,
nontrivial $\ell$-uniform companion lies in $\mathcal B(n,k,\ell)$.
Moreover, if $\mathcal A\subseteq\mathcal B(n,k,\ell)$ and
\[
|\mathcal A|\ge
\binom{n-1}{k-1}-\binom{n-\ell-1}{k-1}+1,
\]
then both $\mathcal A$ and its full companion are nonempty and
nontrivial.

We use the following consequence of
\cite[Theorem~2.8 and Proposition~2.9]{wuzhang26}.

\begin{thm}\label{zhjc}
Let $n>k+\ell$ and $k,\ell\ge2$. Let $\mathcal A\subseteq\binom{[n]}k$
be shifted and satisfy
\[
|\mathcal A|\ge \binom{n-1}{k-1}-\binom{n-\ell-1}{k-1}+1.
\]
If $\overline N_\ell(\mathcal A)$ is nonempty and nontrivial, then
\[
|N_\ell(\mathcal A)|\ge
\left|N_\ell\!\left(\mathcal L_\ell([n],k,|\mathcal A|)\right)\right|.
\]
In particular, equality holds when
$\mathcal A=\mathcal L_\ell([n],k,|\mathcal A|)$.
Moreover, if $\ell=k\ge3$ and $|\mathcal A|=h(n,k,1)$,
equality holds if and only if $\mathcal A$ is a maximum nontrivial
intersecting family.
\end{thm}

\begin{proof}
The neighborhood comparison follows directly from
\cite[Theorem~2.8]{wuzhang26}. We only explain the specialization
of its equality statement to shifted families.

Suppose that $\ell=k\ge3$ and $|\mathcal A|=h=h(n,k,1)$.
Put $\mathcal B=\overline N_k(\mathcal A)$, which is also shifted.
Both families are nontrivial by the preceding observations.
The restricted initial segment of size $h$ is
$\mathcal H(n,k,1)$, whose full companion is itself.
Thus equality is equivalent to $|\mathcal B|=h$, and
\cite[Proposition~2.9]{wuzhang26} applies.

In the first form of that proposition, the distinguished point
has uniquely largest degree in each family, so shiftedness makes
it $1$. Each family's unique member avoiding $1$ must then be
$[2,k+1]$. Hence $\mathcal A=\mathcal B=\mathcal H(n,k,1)$.

The second form occurs only when $k=3$. Write its distinguished
point as $x$ and its two pairs as $P,Q$. If $P\ne Q$, write
$P=\{y,u\}$ and $Q=\{y,v\}$, and choose
$z\notin\{x,y,u,v\}$. Then
$\{y,u,z\},\{x,v,z\}\in\mathcal A$, whereas
$\{x,u,z\},\{y,v,z\}\notin\mathcal A$.
Whichever of $x,y$ comes first, shifting the other to it
contradicts shiftedness. Therefore $P=Q$, and
\[
\mathcal A=\mathcal B
=\{F\in\binom{[n]}3:|F\cap(\{x\}\cup P)|\ge2\},
\]
a maximum nontrivial intersecting family.

Conversely, if $\mathcal A$ is a maximum nontrivial intersecting
family, then $\mathcal A\subseteq\overline N_k(\mathcal A)$.
Theorem~\ref{mh} bounds the latter's size by $h$, so equality holds.
\end{proof}

\section{Weighted bounds for shifted families}

We apply the restricted comparison in Theorem~\ref{zhjc} to prove
the weighted bound when the smaller family is shifted and
nontrivial. We first establish a weighted neighborhood estimate
for families whose members meet a fixed set. Put
\[
\mathcal T_{n,k,p}
=\left\{A\in\binom{[n]}k:A\cap[p]\ne\emptyset\right\}.
\]

\begin{lem}\label{nab}
Let $n,k,p,m$ be positive integers with $n>k+p$ and $m\ge2$. If $\emptyset\ne\mathcal A\subseteq\mathcal T_{n,k,p}$, then
\begin{equation}
\label{nab-bound} |\mathcal A|-m|N_p(\mathcal A)| \le \max\left\{ 1-m\binom{n-k}p,\, \binom nk-\binom{n-p}k-m\left(\binom np-1\right) \right\}.
\end{equation}
Equality holds only if $|\mathcal A|=1$ or $\mathcal A=\mathcal T_{n,k,p}$. Each of these families gives equality precisely when its
corresponding expression attains the displayed maximum.
\end{lem}

\begin{proof}
Since $[p]\in\overline N_p(\mathcal A)$ and
\[
|\mathcal A|-m|N_p(\mathcal A)|
=|\mathcal A|+m|\overline N_p(\mathcal A)|-m\binom np,
\]
it suffices to maximize $|\mathcal X|+m|\mathcal Y|$ over nonempty
cross-intersecting families of uniformities $k$ and $p$.
By Lemma~\ref{lem1}, we may choose a maximizing pair of initial
lexicographic segments. Maximality makes the pair saturated:
each family is the full companion of the other.

We first allow an arbitrary weight $\lambda>0$ and maximize
$|\mathcal X|+\lambda|\mathcal Y|$. We will show that the canonical
indices of $\mathcal X$ are consecutive and that its expansion
has either one term or $p$ terms. Write
\[
|\mathcal X|=\sum_{i=0}^s\binom{n-a_i}{k-a_i+i},
\qquad
1\le a_0<\cdots<a_s,\quad a_i\le k+i,\quad s\le n-k-1.
\]
Here $a_0\ge1$ because the initial segment $\mathcal Y$ contains
$[p]$, so $\mathcal X\subseteq\mathcal T_{n,k,p}$.
Theorem~\ref{cokk} gives
\[
|N_p(\mathcal X)|=\sum_{i=0}^s\binom{n-a_i}{p-i}.
\]
Thus the equivalent objective is
$\rho_\lambda(\mathcal X)=|\mathcal X|-\lambda|N_p(\mathcal X)|$.

First, $s\le p-1$. Otherwise, append the unit
$\binom{n-k-s-1}{n-k-s-1}=1$ to the size expansion.
Its neighborhood contribution is zero because $p-s-1<0$.
If the new lower index is zero, combine the terminal terms
to restore canonical form; the neighborhood sum is unchanged.
The next lexicographic segment therefore has one more member
and the same $p$-neighborhood. It remains in $\mathcal T_{n,k,p}$:
the current segment is proper, since the full domain has a
canonical expansion with exactly $p$ terms.
This contradicts maximality.

Next, suppose that the indices have a gap, and choose the largest
$j\ge1$ with $a_j>a_{j-1}+1$.
Write the consecutive tail as $a_{j+h}=t+h$ for $0\le h<L$,
where $L=s-j+1$. Then $k-t+j\ge0$ and $1\le L\le p-j$.
The tail contributes $A-\lambda B$, where
\[
A=\sum_{h=0}^{L-1}\binom{n-t-h}{k-t+j},
\qquad B=\sum_{h=0}^{L-1}\binom{n-t-h}{p-j-h}.
\]
Deleting the tail leaves a nonempty initial segment.
Alternatively, the gap allows us to decrease all tail indices
by $1$. This changes the objective by $A^+-\lambda B^+$, where
\[
A^+=\sum_{h=0}^{L-1}\binom{n-t-h}{k-t+j+1},
\qquad B^+=\sum_{h=0}^{L-1}\binom{n-t-h}{p-j-h-1}.
\]
Both segments remain in $\mathcal T_{n,k,p}$, since their
expansions have at most $p$ terms and satisfy $a_i\ge i+1$.
We compare the two changes through
\begin{equation}\label{nab-tail-comparison}
\frac{A^+}{A}\ge\frac{n-k-j}{k-t+j+2}
>\frac{p-j}{n-t-p+j+1}\ge\frac{B^+}{B}.
\end{equation}
For the first inequality, the ratios of corresponding summands
in $A^+$ and $A$ decrease with $h$. Extending both sums to
$0\le h\le n-k-j$ therefore gives
\[
\frac{A^+}{A}\ge
\frac{\binom{n-t+1}{k-t+j+2}}{\binom{n-t+1}{k-t+j+1}}
=\frac{n-k-j}{k-t+j+2}.
\]
The middle inequality follows from $n>k+p$.
For the last, the corresponding summand ratios satisfy
\[
\frac{\binom{n-t-h}{p-j-h-1}}{\binom{n-t-h}{p-j-h}}
=\frac{p-j-h}{n-t-p+j+1}
\le\frac{p-j}{n-t-p+j+1}.
\]
Maximality under deletion gives $A-\lambda B\ge0$,
and maximality under the second change gives
$A^+-\lambda B^+\le0$.
Hence $A/B\ge\lambda\ge A^+/B^+$, contradicting
\eqref{nab-tail-comparison}.
The indices are therefore consecutive:
\[
a_i=t+i\qquad(0\le i<L)
\]
for some $1\le t\le k$ and $1\le L\le p$.

For fixed $t$, the objective becomes
\[
\rho_\lambda(\mathcal X)
=\sum_{h=0}^{L-1}
\left(\binom{n-t-h}{k-t}-\lambda\binom{n-t-h}{p-h}\right).
\]
The ratios of the two binomial terms strictly increase with $h$, since
\[
\frac{\binom{n-t-h-1}{k-t}/\binom{n-t-h-1}{p-h-1}}
     {\binom{n-t-h}{k-t}/\binom{n-t-h}{p-h}}
=\frac{n-k-h}{p-h}>1
\qquad(0\le h\le p-2).
\]
Thus, as $L$ increases, successive differences can change sign
only from negative to positive. Every value with $1<L<p$
is strictly below the larger endpoint value, so $L=1$ or $L=p$.
These conclusions hold for every $\lambda>0$.
Interchanging the families and using weight $1/\lambda$ gives
the corresponding conclusions for $\mathcal Y$.

The consecutive expansion describes the family
\[
\mathcal X=
\left\{X\in\binom{[n]}k:
[t-1]\subseteq X,\ X\cap[t,t+L-1]\ne\emptyset\right\}.
\]
Its full companion consists of the $p$-sets meeting $[t-1]$
and the $p$-sets containing $[t,t+L-1]$.
Indeed, any $p$-set satisfying neither condition misses a member
of $\mathcal X$, which can be completed outside that set because
$n>k+p$.

If $L=p\ge2$ and $t\ge2$, the companion consists of all $p$-sets
meeting $[t-1]$, together with $[t,t+p-1]$.
Its canonical indices are
\[
1,2,\ldots,t-1,\ p+t-1,
\]
contradicting the absence of gaps.
Hence $L=p\ge2$ forces $t=1$ and
$\mathcal X=\mathcal T_{n,k,p}$.
The case $p=1$ is already covered by $L=1$.

If $L=1$, then $\mathcal X$ is the full family containing $[t]$,
and $\mathcal Y$ consists of all $p$-sets meeting $[t]$.
The expansion of $\mathcal Y$ has $t$ consecutive terms,
so its endpoint property gives $t=1$ or $t=k$.
The latter gives $|\mathcal X|=1$.
The former gives two full stars through $1$; if $k=1$ or $p=1$,
this is already one of the required endpoints.

It remains to exclude two full stars when $k,p\ge2$, now with
$\lambda=m$. Compare them with the full $k$-family containing
$[2]$ and its companion, and with the full $k$-family meeting
$[2]$ and its companion. The changes in
$|\mathcal X|+m|\mathcal Y|$ are $\delta$ and $-\delta$, where
\[
\delta=m\binom{n-2}{p-1}-\binom{n-2}{k-1}.
\]
If $\delta\ne0$, one pair gives a strict improvement.
If $\delta=0$, both pairs are maximizing.
For $k\ge3$, the first pair has a dual expansion of length $2$,
strictly between $1$ and $k$, which is impossible.
For $k=2$ and $p\ge3$, the second pair has a $k$-side expansion
of length $2$, strictly between $1$ and $p$, again impossible.
Finally, if $k=p=2$, then $\delta=(m-1)(n-2)>0$.
Thus two full stars cannot maximize the objective.

Every maximizing lexicographic pair consequently has
$|\mathcal X|=1$ or $\mathcal X=\mathcal T_{n,k,p}$.
The corresponding objective values are
\[
1+m\left(\binom np-\binom{n-k}p\right),
\qquad \binom nk-\binom{n-p}k+m.
\]
Subtracting $m\binom np$ proves \eqref{nab-bound}.

For an arbitrary equality family, lexicographic replacement
preserves its size and cannot decrease the objective.
Its size must therefore be one of the two endpoint sizes.
In the second case, containment in $\mathcal T_{n,k,p}$
determines the family. The two endpoint counts give the stated
conditions for equality.
\end{proof}

We now apply Lemma~\ref{nab} to the last nonempty minimum-element
slice of a restricted lexicographic family. This gives the
weighted bound needed in the shifted case.

\begin{thm}\label{sab}
Let $n>2k$, $k\ge2$, and $m\ge2$ be integers. Let $\mathcal A,\mathcal B\subseteq\binom{[n]}k$ be cross-intersecting, with
\[
|\mathcal A|\le h(n,k,1)\le|\mathcal B|.
\]
If $\mathcal A$ is shifted and nontrivial, then
\begin{equation}
\label{sab-bound} m|\mathcal A|+|\mathcal B| \le\max\left\{ \sum_{i=2}^{k+1}\binom{k+1}{i}\binom{n-k-1}{k-i}+m(k+1), \ (m+1)h(n,k,1)\right\}.
\end{equation}
Equality holds if and only if at least one of the following holds:
\begin{enumerate}
\item[\textup{(i)}]
\[
\sum_{i=2}^{k+1}\binom{k+1}{i}\binom{n-k-1}{k-i}+m(k+1) \ge(m+1)h(n,k,1),
\]
$\mathcal A=\binom{[k+1]}k$ and $\mathcal B=\mathcal B(n,k,k)$;
\item[\textup{(ii)}]
\[
\sum_{i=2}^{k+1}\binom{k+1}{i}\binom{n-k-1}{k-i}+m(k+1) \le(m+1)h(n,k,1),
\]
and $\mathcal A=\mathcal B$ is a maximum nontrivial intersecting family in $\binom{[n]}k$.
\end{enumerate}
\end{thm}

\begin{proof}
If $k=2$, shiftedness and nontriviality give
$\binom{[3]}2\subseteq\mathcal A$.
Since $h(n,2,1)=3$, we have $\mathcal A=\binom{[3]}2$.
Cross-intersection and $|\mathcal B|\ge3$ then force
$\mathcal B=\mathcal A$.
Both expressions in \eqref{sab-bound} equal $3(m+1)$.
Assume from now on that $k\ge3$.

Choose an admissible pair maximizing $m|\mathcal A|+|\mathcal B|$.
Maximality first gives
\[
\mathcal B=\overline N_k(\mathcal A).
\]
Thus $\mathcal B$ is shifted.
Also, a shifted nontrivial $k$-uniform family contains
$\binom{[k+1]}k$: a member avoiding $1$ can be shifted to
$[2,k+1]$, and further shifts give the other members of the clique.
Applying this to $\mathcal A$ yields
\[
\mathcal B\subseteq\mathcal B(n,k,k).
\]
The family $\mathcal B$ is nontrivial as well.
Otherwise, shiftedness would force every member to contain $1$,
and the preceding containment would give
\[
|\mathcal B|
\le\binom{n-1}{k-1}-\binom{n-k-1}{k-1}=h(n,k,1)-1,
\]
a contradiction.

The full companion of $\mathcal B$ contains $\mathcal A$
and is therefore nontrivial.
Since $|\mathcal B|\ge h(n,k,1)$, Theorem~\ref{mh} bounds that
companion's size by $h(n,k,1)$.
Saturating $\mathcal A$ is therefore admissible, and maximality gives
\[
\mathcal A=\overline N_k(\mathcal B),\qquad
m|\mathcal A|+|\mathcal B|=m\binom nk+\rho(\mathcal B),
\]
where
\[
\rho(\mathcal B)=|\mathcal B|-m|N_k(\mathcal B)|.
\]

The restricted comparison in Theorem~\ref{zhjc} allows us to replace $\mathcal B$ by
$\mathcal L_k([n],k,|\mathcal B|)$ without decreasing $\rho$.
The replacement is admissible: it is shifted, lies in
$\mathcal B(n,k,k)$, and is nontrivial because its size is at
least $h(n,k,1)$.
Its full companion contains $\binom{[k+1]}k$, so it is nontrivial
and has size at most $h(n,k,1)$ by Theorem~\ref{mh}.
We may therefore assume that
\[
\mathcal B=\mathcal L_k([n],k,|\mathcal B|).
\]
This replacement preserves $|\mathcal B|$.
We will use this fact to recover the original equality pairs.

For $1\le i\le k$, define
\[
\mathcal B_i
=\{B\setminus\{i\}:B\in\mathcal B,\ \min B=i\},
\]
and let $r$ be the largest index with $\mathcal B_r\ne\emptyset$.
The first slice has size $h(n,k,1)-1$, so $2\le r\le k$.
Every slice before $r$ is full. Hence
\begin{align}
|\mathcal B|
&=\sum_{i=1}^{r-1}
\left(\binom{n-i}{k-1}-\binom{n-k-1}{k-1}\right)
+|\mathcal B_r|,\label{sab-slices}\\
|N_k(\mathcal B)|
&=\sum_{i=1}^{r-1}
\left(\binom{n-i}{k-i+1}-1\right)
+|N_{k-r+1}(\mathcal B_r)|.\label{sab-neighborhood}
\end{align}
The last neighborhood is taken on $[r+1,n]$.
For \eqref{sab-neighborhood}, partition a neighboring $k$-set
by the first element of $[r]$ that it misses.
If that element is $i<r$, the full $i$th slice gives all
$\binom{n-i}{k-i+1}$ possibilities except $[k+1]\setminus\{i\}$.
If it is $r$, the remaining condition is precisely membership
in $N_{k-r+1}(\mathcal B_r)$.
A set containing $[r]$ has no disjoint member of $\mathcal B$.

Apply Lemma~\ref{nab} to $\mathcal B_r$ on $[r+1,n]$,
with parameters $(n-r,k-1,k-r+1)$.
The required numerical condition follows from $n>2k$.
Both endpoint replacements for this slice are admissible.
The singleton endpoint leaves at least $h(n,k,1)$ members,
and the full endpoint stays within $\mathcal B(n,k,k)$;
their full companions satisfy the required conditions by
Theorem~\ref{mh}.
Together with \eqref{sab-slices}--\eqref{sab-neighborhood},
maximality therefore forces
\[
|\mathcal B_r|=1
\quad\text{or}\quad
|\mathcal B_r|=\binom{n-r}{k-1}-\binom{n-k-1}{k-1}.
\]

Suppose first that $|\mathcal B_r|=1$.
If $r>2$, delete the last slice and call the remaining family
$\mathcal B'$.
Its first two full slices ensure $|\mathcal B'|\ge h(n,k,1)$,
while \eqref{sab-slices}--\eqref{sab-neighborhood} give
\[
\rho(\mathcal B')-\rho(\mathcal B)
=m\binom{n-k-r+1}{k-r+1}-1>0,
\]
contradicting maximality.
Thus $r=2$ and
\[
\mathcal B=
\left\{B\in\binom{[n]}k:
1\in B,\ B\cap[2,k+1]\ne\emptyset\right\}
\cup\{[2,k+1]\}.
\]
Its full companion is itself, so
$|\mathcal A|=|\mathcal B|=h(n,k,1)$.

Now suppose that the last slice is full.
If $2<r<k$, compare $\mathcal B$ with the families consisting
of the first $r-1$ and the first $r+1$ full slices.
Both are admissible, so maximality gives
\begin{align*}
\binom{n-r}{k-1}-\binom{n-k-1}{k-1}
&\ge m\left(\binom{n-r}{k-r+1}-1\right),\\
\binom{n-r-1}{k-1}-\binom{n-k-1}{k-1}
&\le m\left(\binom{n-r-1}{k-r}-1\right).
\end{align*}
All denominators below are positive. It follows that
\begin{equation}\label{sab-forced-ratio}
\frac{\binom{n-r-1}{k-r}-1}
     {\binom{n-r-1}{k-1}-\binom{n-k-1}{k-1}}
\ge
\frac{\binom{n-r}{k-r+1}-1}
     {\binom{n-r}{k-1}-\binom{n-k-1}{k-1}}.
\end{equation}
On the other hand,
\[
\frac{\binom{n-r-1}{k-r}-1}
     {\binom{n-r-1}{k-1}-\binom{n-k-1}{k-1}}
=\frac{\displaystyle\sum_{i=n-k-1}^{n-r-2}\binom{i}{n-k-2}}
       {\displaystyle\sum_{i=n-k-1}^{n-r-2}\binom{i}{k-2}}
<\frac{\displaystyle\sum_{i=n-k-1}^{n-r-1}\binom{i}{n-k-2}}
       {\displaystyle\sum_{i=n-k-1}^{n-r-1}\binom{i}{k-2}}
=\frac{\binom{n-r}{k-r+1}-1}
        {\binom{n-r}{k-1}-\binom{n-k-1}{k-1}},
\]
contradicting \eqref{sab-forced-ratio}.
The strict inequality holds because
$\binom{i}{n-k-2}/\binom{i}{k-2}$ strictly increases with $i$
when $n>2k$.

If $r=2$, retain only the first full slice and the first member
of the second slice, obtaining the Hilton--Milner family
$\mathcal B'$ above. Then
\[
\rho(\mathcal B)-\rho(\mathcal B')
=-(m-1)\left(\binom{n-2}{k-1}
-\binom{n-k-1}{k-1}-1\right)<0,
\]
where the strictness follows from $k\ge3$ and $n>2k$.
This again contradicts maximality.
The full-slice case therefore has $r=k$, giving
\[
\mathcal B=\mathcal B(n,k,k),\qquad
\mathcal A=\binom{[k+1]}k.
\]
The two resulting configurations are admissible, and their sizes
give the two values in \eqref{sab-bound}.

Finally, let the original pair attain equality.
As above, saturation makes $\mathcal B$ shifted and places it
in $\mathcal B(n,k,k)$.
Lexicographic replacement preserves its size, so the preceding
argument gives $|\mathcal B|=h(n,k,1)$ or
$|\mathcal B|=|\mathcal B(n,k,k)|$.

In the first case, $|\mathcal A|\le h(n,k,1)$ and equality in
\eqref{sab-bound} force $|\mathcal A|=h(n,k,1)$ and the
comparison in \textup{(ii)}. Hence
\[
|N_k(\mathcal B)|=\binom nk-h(n,k,1)
=|N_k(\mathcal H(n,k,1))|.
\]
Since $\mathcal L_k([n],k,h(n,k,1))=\mathcal H(n,k,1)$,
the shifted equality conclusion of Theorem~\ref{zhjc} makes the original
$\mathcal B$ a maximum nontrivial intersecting family.
Thus $\mathcal B\subseteq\overline N_k(\mathcal B)=\mathcal A$,
and equal sizes give $\mathcal A=\mathcal B$.

In the second case, containment forces
$\mathcal B=\mathcal B(n,k,k)$.
Saturation then gives $\mathcal A=\binom{[k+1]}k$ and the
comparison in \textup{(i)}.
Conversely, each listed configuration attains the bound when
its corresponding value is maximal.
\end{proof}

\section{Critically shifted families}

We first prove a weighted bound for two families of different
uniformities. We then combine it with local generator contractions
to handle the critically shifted case of Theorem~\ref{mth}.

\subsection{The two-family bound}

We begin with the shifted case. Nontriviality forces a clique,
which is preserved by the exchanges used below.

\begin{lem}\label{sab1-clique}
Let $n,k,\ell$ be positive integers with $n>k+\ell>2k$, and let $m>0$ be real. Let
$\mathcal{A}\subseteq\binom{[n]}{k}$ and $\mathcal{B}\subseteq\binom{[n]}{\ell}$ be cross-intersecting. Suppose that $\mathcal{A}$ is
nonempty, nontrivial, and shifted. Then
\[
m|\mathcal{A}|+|\mathcal{B}| \leq\max\left\{ m(k+1)+\sum_{i=2}^{k+1}\binom{k+1}{i}\binom{n-k-1}{\ell-i},\ m\binom nk\right\}.
\]
If $m\binom nk>m(k+1)+\sum_{i=2}^{k+1}\binom{k+1}{i}\binom{n-k-1}{\ell-i}$, equality holds if and only if $\mathcal{A}=\binom{[n]}{k}$ and
$\mathcal{B}=\emptyset$.
\end{lem}

\begin{proof}
Put $S=[k+1]$ and $N=n-k-1$.
Shiftedness and nontriviality imply $\binom Sk\subseteq\mathcal A$:
a member avoiding $1$ shifts to $[2,k+1]$, whose shifts give the whole clique.
It suffices to maximize
\[
W=m|\mathcal X|+|\mathcal Y|
\]
over cross-intersecting pairs with
$\binom Sk\subseteq\mathcal X\subseteq\binom{[n]}k$.
Allowing $\mathcal X$ to be unshifted can only increase the maximum.

Among maximizing pairs, choose one with $|\mathcal X|$ largest.
Both families are saturated.
For any permutation $\pi$ preserving $S$ setwise, the pairs
\begin{equation}\label{generator-pair-exchanges}
(\mathcal X\cup\pi\mathcal X,\mathcal Y\cap\pi\mathcal Y),\qquad
(\mathcal X\cap\pi\mathcal X,\mathcal Y\cup\pi\mathcal Y)
\end{equation}
are cross-intersecting, and their first families contain $\binom Sk$.
Both pairs are therefore admissible.
Their objectives add to $2W$ and neither exceeds $W$, so both
are maximizing.
The choice of $|\mathcal X|$ forces $\mathcal X=\pi\mathcal X$.
Saturation then also gives $\mathcal Y=\pi\mathcal Y$.

Now shift both families simultaneously.
The clique constraint and both sizes are preserved, so the resulting
pair has the same invariance.
Membership in $\mathcal X$ therefore depends only on $|X\cap S|$.
Shifting a point into $S$ shows that the allowed intersection sizes
form an upper interval. Hence, for some $0\leq r\leq k$,
\[
\mathcal X=\left\langle\binom Sr\right\rangle_k
 =\{X\in\tbinom{[n]}k:|X\cap S|\geq r\}.
\]
By Lemma~\ref{generator-rules}\textup{(i)}, its full companion is
\[
\mathcal Y=\{Y\in\tbinom{[n]}\ell:|Y\cap S|\geq k+2-r\}.
\]
This includes $\mathcal Y=\emptyset$ for $r=0$.
The corresponding objective is
\[
W_r=m\sum_{i=r}^{k}\binom{k+1}i\binom N{k-i}
       +\sum_{i=k+2-r}^{k+1}\binom{k+1}i\binom N{\ell-i}.
\]
For $0\leq r<k$,
\begin{equation}\label{sab1-generator-difference}
W_{r+1}-W_r
 =\binom{k+1}r\binom N{k-r}(R_r-m),\qquad
R_r=\frac{\binom N{\ell-k-1+r}}{\binom N{k-r}}.
\end{equation}
All binomial coefficients here are positive because $N\geq\ell>k$.
For $0\leq r<k-1$,
\[
\frac{R_{r+1}}{R_r}
 =\frac{(N-\ell+k+1-r)(N-k+r+1)}{(\ell-k+r)(k-r)}>1.
\]
Indeed, the first numerator factor exceeds $k-r$, and the second exceeds $\ell-k+r$.
Thus the differences in \eqref{sab1-generator-difference} can change sign only from negative to positive.
Every interior value is strictly below the larger endpoint value. The endpoints are
\[
W_0=m\binom nk,\qquad
W_k=m(k+1)+\sum_{i=2}^{k+1}\binom{k+1}i\binom{n-k-1}{\ell-i}.
\]
This proves the asserted inequality, also when $k=1$, for which there are no interior values.

Suppose that $W_0>W_k$.
To recover equality for arbitrary pairs, choose a maximizing pair
with $|\mathcal Y|$ largest.
Using the second exchange in \eqref{generator-pair-exchanges} shows that $\mathcal Y$ is invariant
under permutations preserving $S$, and saturation gives the same for $\mathcal X$.
Shifting preserves this secondary maximum. The preceding reduction again gives a threshold pair,
which must have $r=0$ by the strict endpoint comparison.
Hence even the largest possible $|\mathcal Y|$ among maximizing pairs is zero.
Every maximizing pair consequently has $\mathcal Y=\emptyset$ and
$\mathcal X=\binom{[n]}k$. Conversely, this pair attains $W_0$.
\end{proof}

\begin{thm}\label{m2n}
Let $n,m,k,\ell$ be positive integers satisfying $m\geq2$, $n>k+\ell>2k$, and $m\binom{n}{k}<\binom{n}{\ell}$. Let
$\mathcal{A}\subseteq\binom{[n]}{k}$ and $\mathcal{B}\subseteq\binom{[n]}{\ell}$ be cross-intersecting. Suppose $\mathcal{A}$ is nonempty
and nontrivial. Then
\[
m|\mathcal{A}|+|\mathcal{B}| \leq\max\left\{\binom{n}{\ell}-2\binom{n-k}{\ell} +\binom{n-2k}{\ell}+2m,m\binom{n}{k}\right\}.
\]
Equality holds if and only if one of the following occurs:
\begin{enumerate}
\item[(i)] $\binom{n}{\ell}-2\binom{n-k}{\ell}
+\binom{n-2k}{\ell}+2m\geq m\binom{n}{k}$, $\mathcal{A}=\{A_1,A_2\}$ for two disjoint members $A_1,A_2\in\binom{[n]}{k}$, and
\[
\mathcal{B}=\left\{B\in\binom{[n]}{\ell}: B\cap A_1\neq\emptyset,\ B\cap A_2\neq\emptyset\right\}.
\]
\item[(ii)] $\binom{n}{\ell}-2\binom{n-k}{\ell}
+\binom{n-2k}{\ell}+2m\leq m\binom{n}{k}$, $\mathcal{A}=\binom{[n]}{k}$, and $\mathcal{B}=\emptyset$.
\end{enumerate}
\end{thm}
\begin{proof}
Put
\[
M=\binom{n}{\ell}-2\binom{n-k}{\ell} +\binom{n-2k}{\ell}+2m.
\]
We induct on $k$.
Call a pair \emph{extremal} if it maximizes
$m|\mathcal A|+|\mathcal B|$ under the hypotheses.
Every extremal pair is saturated, since adding members preserves
the nontriviality of $\mathcal A$.
The two constructions in the statement have values $M$ and
$m\binom nk$, so the extremal value is at least their maximum.

Suppose first that $k=1$, and write $|\mathcal A|=s\geq2$. By saturation,
\[
|\mathcal B|=\binom{n-s}{\ell-s}\quad\text{if }s\leq\ell, \qquad \mathcal B=\emptyset\quad\text{if }s>\ell.
\]
Let $f(s)=ms+\binom{n-s}{\ell-s}$, where the binomial term is understood to be zero for $s>\ell$. For $2\leq s\leq\ell-1$,
\[
f(s+2)-2f(s+1)+f(s)=\binom{n-s-2}{\ell-s}>0.
\]
Thus $f$ is strictly convex on the integer interval $[2,\ell+1]$, and $f(s)=ms$ is strictly increasing for $s\geq\ell+1$. Since $n>\ell+1$,
a maximum can occur only at $s=2$ or $s=n$. These give (i) and (ii), respectively.

Assume that $k\geq2$ and the result holds for smaller uniformities.
If $\mathcal A$ is shifted, Lemma~\ref{sab1-clique} gives
\[
m|\mathcal A|+|\mathcal B|
\leq\max\left\{
m(k+1)+\sum_{i=2}^{k+1}\binom{k+1}{i}\binom{n-k-1}{\ell-i},
\ m\binom nk\right\}.
\]
Lemma~\ref{inma} makes the first expression strictly smaller than
$M$, proving the bound in this case.
An extremal shifted pair must therefore attain $m\binom nk$.
The equality statement in Lemma~\ref{sab1-clique} then gives
\textup{(ii)}.

Now take an extremal pair. If $\mathcal B=\emptyset$, saturation
gives \textup{(ii)}.
Otherwise, apply simultaneous shifts that change $\mathcal A$
and preserve its nontriviality until no such shift remains.
The process terminates because every effective shift decreases
the sum of the elements in the members of $\mathcal A$.
Write $\mathcal A^\ast$ for the resulting family.
If it is shifted, the preceding argument applies.
Otherwise, every shift that changes it makes it trivial.
Choose such a shift, relabel its indices as $1,2$, and continue
to write $\mathcal A,\mathcal B$ for the extremal pair.
Thus $s_{12}(\mathcal A)$ is trivial.

For $P\subseteq[2]$ and $\mathcal F\subseteq2^{[n]}$, use the section notation
\[
\mathcal F(P)=\{F\setminus P:F\in\mathcal F,\ F\cap[2]=P\}.
\]
We write $\mathcal F(1)$, $\mathcal F(2)$, and $\mathcal F(12)$ for the sections indexed by $\{1\}$, $\{2\}$, and $\{1,2\}$, respectively.
Since $s_{12}$ creates a common element, we have
\begin{equation}
\label{m2n-sections} \mathcal A(\emptyset)=\emptyset,\qquad \mathcal A(1)\cap\mathcal A(2)=\emptyset,\qquad \mathcal A(1),\mathcal A(2)\neq\emptyset.
\end{equation}
Saturation also gives $\mathcal B(12)=\binom{[3,n]}{\ell-2}$.

Next apply simultaneous shifts $s_{1j},s_{2j}$ with $j\in[3,n]$
as long as they change $\mathcal A$ and preserve its nontriviality.
They only move members from a single-index section to the double
section, so \eqref{m2n-sections} remains valid.
In particular, a single-index section cannot become empty without
making $\mathcal A$ trivial.
The resulting pairs remain extremal and saturated. Put
\[
C=\bigcap_{E\in\mathcal A(1)}E,\qquad D=\bigcap_{E\in\mathcal A(2)}E,\qquad \mathcal H=\mathcal A(12).
\]
If $j\in [3,n]\setminus C$, some member of $\mathcal A(1)$ avoids $j$. The corresponding member of $\mathcal A$ avoids both $2$ and $j$ and
is fixed by $s_{2j}$. As a shift can create a common element only at its target index, $s_{2j}(\mathcal A)$ is nontrivial and hence equals
$\mathcal A$. Consequently,
\begin{equation}
\label{m2n-deletion} E\setminus\{j\}\in\mathcal H \qquad(E\in\mathcal A(1),\ j\in E\setminus C).
\end{equation}
The analogous statement holds for $\mathcal A(2)$ and $D$.

If $C=\emptyset$, we can reduce to the shifted case.
Let $\partial$ denote the shadow of codimension one.
Equation~\eqref{m2n-deletion} gives
$\partial\mathcal A(1)\subseteq\mathcal H$.
Shifts on $[3,n]$ preserve this inclusion because
\[
\partial s_{ij}(\mathcal F) \subseteq s_{ij}(\partial\mathcal F).
\]
It also implies $\bigcap_{E\in\mathcal H}E=\emptyset$, because $\mathcal H$ contains the full shadow of a $(k-1)$-set. Thus shifts on
$[3,n]$ preserve the nontriviality of $\mathcal A$. Apply such shifts until its sections are shifted on $[3,n]$. Both single-index sections
then contain the first $(k-1)$-subset $E_0$ of $[3,n]$.

Now apply shifts $s_{1j}$ until they fix $\mathcal A$. The member $\{2\}\cup E_0$ cannot move, since $E_0\setminus\{j\}\in\mathcal H$ for
every $j\in E_0$. The first single-index section is unchanged, while $\mathcal H$ can only increase. Hence these shifts preserve
nontriviality, and at their termination
$\partial\mathcal A(2)\subseteq\mathcal H$ as well.

Apply $s_{12}$. Its new single-index sections
are $\mathcal A(1)\cup\mathcal A(2)$ and $\mathcal A(1)\cap\mathcal A(2)$, both containing $E_0$; the shadow of each is contained in
$\mathcal H$.

Finally, shift on $[3,n]$ again until the sections are shifted there. These shadow inclusions and the inclusion of the second
single-index section in the first are preserved. Every $s_{1j}$ and $s_{2j}$ now fixes $\mathcal A$, as does $s_{12}$. Thus $\mathcal A$ is
shifted and remains nontrivial, so the shifted case applies. The case $D=\emptyset$ follows by interchanging $1$ and $2$. We may therefore
assume $C,D\neq\emptyset$.

\medskip
\noindent\textit{Claim.} Each single-index section is generated by a single $(k-1)$-set.

It is enough to consider $\mathcal A(1)$.
Suppose that it has more than one member, so $1\leq c=|C|\leq k-2$.
Choose $E\in\mathcal A(1)$ and
$G\in\mathcal G_{k-1}(\mathcal A(1))$ with $G\subseteq E$,
on the ground set $[3,n]$.
Delete the other generators and enlarge $G$ to $E$.
This contracts the section to $\{E\}$.
Let $\mathcal A'$ be the resulting family, and obtain $\mathcal B'$
by replacing $\mathcal B(2)$ with all $(\ell-1)$-sets in $[3,n]$
meeting $E$. Keep every other section fixed.

This generator contraction preserves the required shadow constraint. Indeed,
\[
\{E\setminus\{j\}:j\in E\setminus C\}\subseteq\mathcal H
\]
by \eqref{m2n-deletion}. Every new member in the $2$-section of $\mathcal B'$
meets $\{1\}\cup E$ and shares $2$ with the $2$- and $12$-sections of $\mathcal A'$.
All other cross-intersections are preserved because $\mathcal A'\subseteq\mathcal A$.
Also, $E\setminus C\neq\emptyset$, and deletion gives
\[
\bigcap_{H\in\mathcal H}H
\subseteq\bigcap_{j\in E\setminus C}(E\setminus\{j\})=C.
\]
Put $K=\bigcap\mathcal H$. Original nontriviality gives $C\cap D\cap K=\emptyset$.
Since $K\subseteq C\subseteq E$, it follows that
\[
\bigcap\mathcal A'=E\cap D\cap K=D\cap K=C\cap D\cap K=\emptyset.
\]
Thus $\mathcal A'$ remains nontrivial and the replacement is admissible.

We now compare the loss in $\mathcal A(1)$ with the new members
allowed in $\mathcal B(2)$.
Put
\[
N=n-c-2,\qquad p=k-c-1,\qquad q=\ell-1.
\]
For this count, identify $[3,n]\setminus C$ with $[N]$ so that $E\setminus C=[p]$,
and write $\mathcal R=\{F\setminus C:F\in\mathcal A(1)\}$ in these coordinates.
For $1\leq i\leq p$, put $V_i=[i+1,N]$ and
\[
\mathcal R_i=
\{R\setminus[i-1]:R\in\mathcal R\setminus\{[p]\},\ \min([p]\setminus R)=i\}.
\]
Thus $\mathcal R_i\subseteq\binom{V_i}{p-i+1}$ and
$\sum_{i=1}^{p}|\mathcal R_i|=|\mathcal A(1)|-1$.
Define
\[
\mathcal Q_i=
\left\{\{i\}\cup Q:Q\in\binom{V_i}{q-1},\ Q\cap R=\emptyset
\text{ for some }R\in\mathcal R_i\right\}.
\]
If $R\in\mathcal R_i$, the corresponding old section member is $C\cup[i-1]\cup R$.
For $Q\cap R=\emptyset$ as above, the set $\{i\}\cup Q$ is disjoint from this old member:
it avoids $C\cup[i-1]$, and $i\notin R$. It also meets $E$, since $i\in E\setminus C$.
Consequently every member of $\mathcal Q_i$ belongs to $\mathcal B'(2)\setminus\mathcal B(2)$.
The families $\mathcal Q_i$ are pairwise disjoint, since their members have distinct least coordinates.

Count disjoint pairs between $\mathcal R_i$ and $\binom{V_i}{q-1}$.
Each member of $\mathcal R_i$ has $\binom{N-p-1}{q-1}$ such partners, whereas a partner
has at most $\binom{N-i-q+1}{p-i+1}$ neighbors in $\mathcal R_i$. Hence
\[
|\mathcal Q_i|
\geq\frac{\binom{N-p-1}{q-1}}{\binom{N-i-q+1}{p-i+1}}|\mathcal R_i|
=\frac{\binom{N-i}{q-1}}{\binom{N-i}{p-i+1}}|\mathcal R_i|.
\]
The last ratio increases with $i$: its successive ratio is
$(N-i-q+1)/(p-i+1)>1$, since $N>p+q$.
It follows that
\begin{equation}\label{m2n-generator-gain}
|\mathcal B'(2)|-|\mathcal B(2)|
\geq\Lambda_c\bigl(|\mathcal A(1)|-1\bigr),\qquad
\Lambda_c=\frac{\binom{n-c-3}{\ell-2}}{\binom{n-c-3}{k-c-1}}.
\end{equation}

To obtain a strict improvement, it remains to show $\Lambda_c>m$.
The quantities $\Lambda_c$ increase with $c$, since
\[
\frac{\Lambda_{c+1}}{\Lambda_c}
=\frac{n-\ell-c-1}{k-c-1}>1\qquad(1\leq c<k-2).
\]
At $c=1$,
\[
\frac{\Lambda_1}{\binom n\ell/\binom nk}
=\frac{\ell(\ell-1)(n-\ell)(n-\ell-1)}
       {k(k-1)(n-k)(n-k-1)}>1.
\]
For the strict inequality, both
$\ell(n-\ell)>k(n-k)$ and
$(\ell-1)(n-\ell-1)>(k-1)(n-k-1)$ follow from $\ell>k$ and $n>k+\ell$.
Thus $\Lambda_c\geq\Lambda_1>\binom n\ell/\binom nk>m$.
Since the only changed sections are $\mathcal A(1)$ and $\mathcal B(2)$,
\eqref{m2n-generator-gain} yields
\[
m|\mathcal A'|+|\mathcal B'|-m|\mathcal A|-|\mathcal B|
\geq(\Lambda_c-m)\bigl(|\mathcal A(1)|-1\bigr)>0,
\]
contrary to extremality.
Thus $\mathcal A(1)=\{C\}$ with $|C|=k-1$.
The same argument gives $\mathcal A(2)=\{D\}$ and $|D|=k-1$.
This proves the claim.

Saturation now gives
\[
|\mathcal B(1)|=|\mathcal B(2)| =\binom{n-2}{\ell-1}-\binom{n-k-1}{\ell-1}.
\]

If $\mathcal H=\emptyset$, nontriviality gives $C\cap D=\emptyset$,
so $\mathcal A$ consists of two disjoint $k$-sets.
Saturation then gives \textup{(i)}.
It remains to exclude $\mathcal H\ne\emptyset$.
We will prove the stronger bound
\begin{equation}
\label{aim}
m|\mathcal A|+|\mathcal B|
\leq\binom{n}{\ell}-2\binom{n-k}{\ell}
-\binom{n-k-1}{\ell-1}+\binom{n-2k+1}{\ell}+3m,
\end{equation}
with equality if and only if $|C\cap D|=1$ and $\mathcal H=\{C\setminus D\}$ or $\mathcal H=\{D\setminus C\}$.

For the rest of this proof, neighborhoods are taken on $[3,n]$.
Put $v=n-k$ and $\alpha=|C\setminus D|+1$, and define
\[
J=|N_\ell(\mathcal H)\cap N_\ell(\{C,D\})|.
\]
Since $\mathcal B$ is saturated,
\[
|\mathcal B(\emptyset)| =\binom{n-2}{\ell}-2\binom{n-k-1}{\ell} +\binom{v-\alpha}{\ell}-|N_\ell(\mathcal H)|+J.
\]
Combining this with the counts for the other sections gives
\[
|\mathcal B| =\binom{n}{\ell}-2\binom{n-k}{\ell} +\binom{v-\alpha}{\ell}-|N_\ell(\mathcal H)|+J.
\]
Define
\[
\Delta=|N_\ell(\mathcal H)|-m|\mathcal H|-J-\binom{v-\alpha}{\ell}+\binom{v-k}{\ell},
\]
and
\[
\Gamma=\binom{v-1}{\ell-1}-\binom{v-k}{\ell-1}.
\]
Since $|\mathcal A|=2+|\mathcal H|$, we have
\begin{equation}
\label{m2n-objective-gap} m|\mathcal A|+|\mathcal B|=M-\Delta.
\end{equation}
In this notation, \eqref{aim} is equivalent to
\begin{equation}
\label{m2n-final-gap} \Delta\geq\Gamma-m.
\end{equation}

Let $K=\bigcap_{E\in\mathcal H}E$, and put
\[
x=|K\cap(C\setminus D)|+1,\qquad y=|K\cap(D\setminus C)|+1,\qquad z=|K\setminus(C\cup D)|,\qquad t=|K|+2.
\]
Nontriviality gives $K\cap C\cap D=\emptyset$. Hence
\[
1\leq x,y\leq \alpha,\qquad x+y+z=t,\qquad 2\leq t\leq k.
\]
Also, $C\neq D$ by \eqref{m2n-sections}, so $\alpha\geq2$. Every set counted by $J$ avoids $K$ and at least one of $C,D$. Therefore
\begin{equation}
\label{m2n-J} J\leq\binom{v-x-z}{\ell}+\binom{v-y-z}{\ell} -\binom{v-\alpha-z}{\ell}.
\end{equation}
Since $x\leq \alpha$, we have
\[
\binom{v-\alpha}{\ell}-\binom{v-\alpha-z}{\ell} \leq\binom{v-x}{\ell}-\binom{v-x-z}{\ell}.
\]
Moreover, discrete convexity and $x,y+z\geq1$ give
\[
\binom{v-x}{\ell}+\binom{v-y-z}{\ell} \leq\binom{v-1}{\ell}+\binom{v-t+1}{\ell}.
\]
Combining the two estimates with \eqref{m2n-J} gives
\begin{equation}
\label{m2n-J-combined} J+\binom{v-\alpha}{\ell} \leq\binom{v-1}{\ell}+\binom{v-t+1}{\ell}.
\end{equation}

First suppose that $|\mathcal H|\geq2$. Then $2\leq t\leq k-1$. Put $q=k-t\geq1$ and
\[
\mathcal H'=\{E\setminus K:E\in\mathcal H\} \subseteq\binom{[3,n]\setminus K}{q}.
\]
This family is nonempty and nontrivial, and $|[3,n]\setminus K|=n-t$. Furthermore,
\[
\frac{\binom{n-t}{\ell}}{\binom{n-t}{q}} =\frac{\binom{n}{\ell}}{\binom{n}{k}} \prod_{j=0}^{t-1}\frac{n-\ell-j}{k-j} >\frac{\binom{n}{\ell}}{\binom{n}{k}}>m.
\]
Since $n-t>q+\ell>2q$, induction applies to $\mathcal H'$ and the family of all $\ell$-subsets of $[3,n]\setminus K$ meeting every member of
$\mathcal H'$. The $\ell$-neighborhood of $\mathcal H'$ in $[3,n]\setminus K$ is exactly $N_\ell(\mathcal H)$, so
\[
|N_\ell(\mathcal H)|-m|\mathcal H|
\geq\min\left\{2\binom{v}{\ell}-\binom{v-q}{\ell}-2m,\,
\binom{n-t}{\ell}-m\binom{n-t}{q}\right\}.
\]
Lemma~\ref{m2n-binomial}\textup{(ii)} makes the first term
strictly smaller than the second. Thus
\[
|N_\ell(\mathcal H)|-m|\mathcal H| \geq2\binom{v}{\ell}-\binom{v-q}{\ell}-2m.
\]
Since $q+t=k$,
\[
\binom{v}{\ell}-\binom{v-q}{\ell} -\binom{v-t+1}{\ell}+\binom{v-k+1}{\ell} =\sum_{i=1}^{q}\sum_{j=1}^{t-1} \binom{v-i-j}{\ell-2}.
\]
Together with \eqref{m2n-J-combined}, this yields
\[
\Delta\geq\Gamma+\sum_{i=1}^{q}\sum_{j=1}^{t-1}
\binom{v-i-j}{\ell-2}-2m
\geq\Gamma+\binom{v-2}{\ell-2}-2m>\Gamma-m,
\]
where the last inequality follows from Lemma~\ref{m2n-binomial}(i).

Now suppose that $\mathcal H=\{K\}$. Then $t=k$ and $|N_\ell(\mathcal H)|=\binom{v}{\ell}$. By
\eqref{m2n-J-combined},
\[
\Delta\geq\binom{v}{\ell}-m-\binom{v-1}{\ell} -\binom{v-k+1}{\ell}+\binom{v-k}{\ell} =\Gamma-m.
\]
This proves \eqref{m2n-final-gap} and hence \eqref{aim}.

Equality in \eqref{aim} requires $\mathcal H=\{K\}$.
Then \eqref{m2n-J} is exact.
Since $\alpha\geq2$ and $\binom{v-3}{\ell-2}>0$,
equality in \eqref{m2n-J-combined} holds precisely when
\[
z=0\quad\text{and}\quad(x=1\text{ or }y=1).
\]
Indeed, the convexity inequality is strict when $x>1$ and $y+z>1$. If $x=1$ and $z>0$, the preceding binomial-difference inequality is
strict. Thus equality is equivalent to
\begin{equation}
\label{m2n-singleton-equality} K\subseteq C\setminus D\quad\text{or}\quad K\subseteq D\setminus C.
\end{equation}

We also need $C\cap D\ne\emptyset$.
In this equality case, $|\mathcal A|=3$, so the family
$\mathcal A^\ast$ obtained after the initial shifts has three members.
It must be intersecting.
Otherwise, write $\mathcal A^\ast=\{F_1,F_2,F_3\}$ with
$F_1\cap F_2=\emptyset$, and choose
$p\in F_1\setminus F_3$ and $q\in F_2\setminus F_3$.
If $p<q$, the shift $s_{pq}$ moves $F_2$ to
$(F_2\setminus\{q\})\cup\{p\}$ and fixes the other members.
It is effective because $k\geq2$ and $p\notin F_3$,
and the resulting common intersection remains empty.
If $q<p$, use $s_{qp}$ instead.
Either shift contradicts the stopping property of $\mathcal A^\ast$.
Relabeling and subsequent shifts preserve intersection, so the
current $\mathcal A$ is intersecting. Hence $C\cap D\ne\emptyset$.

By symmetry, assume in \eqref{m2n-singleton-equality} that $K\subseteq C\setminus D$. Since $|K|=k-2$ and $|C|=|D|=k-1$, we have $|C\cap
D|\leq1$. Consequently, $|C\cap D|=1$ and $K=C\setminus D$. The other case gives $K=D\setminus C$. Conversely, either configuration gives
equality by direct substitution. This proves the equality characterization in
\eqref{aim}.

Finally, Lemma~\ref{m2n-binomial}(i) gives
\[
\Gamma=\sum_{i=2}^{k}\binom{v-i}{\ell-2} \geq\binom{v-2}{\ell-2}>m.
\]
Hence \eqref{m2n-final-gap} implies $\Delta>0$, and
\eqref{m2n-objective-gap} yields
$m|\mathcal A|+|\mathcal B|<M$ whenever $\mathcal H\neq\emptyset$. This contradicts extremality and completes the induction.

To recover the original equality pairs, note that all shifts and
relabelings preserve both sizes.
If the resulting pair has form \textup{(i)}, the original
$\mathcal A$ also has two members; nontriviality makes them disjoint,
and saturation determines $\mathcal B$.
If it has form \textup{(ii)}, the original sizes already force
$\mathcal A=\binom{[n]}k$ and $\mathcal B=\emptyset$.
Conversely, the two constructions attain equality when their
respective values are maximal.
\end{proof}

\subsection{The main theorem}

We first describe the local structure at a critical pair.
The order on the remaining coordinates will be used for the
generator contraction in the main proof.

\begin{lem}\label{critical-local-generators}
Let $\mathcal A\subseteq\binom{[n]}k$ be nontrivial and critically shifted.
There is a critical pair containing $1$. Choose the least $j$ for which $s_{1j}(\mathcal A)$ is trivial,
and then name these two coordinates $1,2$, retaining the induced order on the remaining coordinates.
Put
\[
C=\bigcap\mathcal A(1),\qquad D=\bigcap\mathcal A(2),\qquad
\mathcal H=\mathcal A(12).
\]
Then $\mathcal A(\emptyset)=\emptyset$, the two single-index sections are nonempty and disjoint,
and the following hold:
\begin{enumerate}
\item[\textup{(i)}] $C,D\neq\emptyset$ and $|C\cup D|\geq k$.
\item[\textup{(ii)}] If $E\in\mathcal A(1)$ and $x\in E\setminus C$, then $E\setminus\{x\}\in\mathcal H$;
the analogous assertion holds for $\mathcal A(2),D$.
\item[\textup{(iii)}] The original family is fixed by every left shift within the remaining coordinates outside $C$,
and likewise outside $D$. Consequently
\[
\{E\setminus C:E\in\mathcal A(1)\},\qquad
\{E\setminus D:E\in\mathcal A(2)\}
\]
are shifted on their respective remaining ground sets.
Moreover, $D\setminus C$ and $C\setminus D$ are initial intervals in those respective orders.
\item[\textup{(iv)}] Every $H\in\mathcal H$ contains $C\setminus D$ or $D\setminus C$.
\end{enumerate}
\end{lem}

\begin{proof}
For a nontrivial family, a shift along $\{u,v\}$ creates a common point precisely when every member
meets $\{u,v\}$ and its two single-index sections are disjoint. This criterion is independent of the orientation.
In particular, a pair avoided by a member cannot be critical, and its left shift fixes a critically shifted family.

Work first in the original coordinate order.
If no critical pair contained $1$, every shift with target $1$ would fix $\mathcal A$.
Choose a critical pair $\{u,v\}$ and a member avoiding $1$.
It must contain both $u,v$, since otherwise a fixed shift to $1$ would produce a member avoiding the pair.
Replacing $u$ and $v$ separately by $1$ then gives a common member in the two sections at $\{u,v\}$,
a contradiction. Thus the required $j$ exists.
Write $C,D$ temporarily for the cores of the sections at $1,j$.

We first show $[2,j-1]\subseteq C$.
If $x<j$ and $x\notin C$, a member of the $1$-section avoids $x,j$, so $s_{xj}$ fixes $\mathcal A$.
Every member of the $j$-section must then contain $x$, since otherwise that fixed shift would remove
its only point of the critical pair. Minimality of $j$ also makes $s_{1x}$ fixed.
Starting with $\{j\}\cup E$ in the $j$-section, apply $s_{1x}$ and then $s_{xj}$.
The result $\{1\}\cup E$ lies in the other section, contradicting their disjointness.

If $x\notin C$ is a remaining coordinate, then $x>j$.
A member of the $1$-section avoiding $x$ shows that $s_{jx}$ is fixed, which proves the first deletion property.
The second follows from the fixed shift $s_{1x}$ whenever $x\notin D$.
Now take two remaining coordinates $u<v$ outside $C$ and choose $E\in\mathcal A(1)$ omitting $u$.
If $v\notin E$, the member $\{1\}\cup E$ avoids the pair.
If $v\in E$, the deletion property supplies $\{1,j\}\cup(E\setminus\{v\})$, again avoiding the pair.
Thus $s_{uv}$ is fixed. The same reasoning applies outside $D$.
Now rename $j$ as $2$, retaining the order of the other coordinates.
The local shifting and deletion properties just proved remain valid
in this inherited order.

If $C$ or $D$ were empty, both single-index sections would be shifted on the whole remaining ground set,
so both would contain its first $(k-1)$-set. Thus both cores are nonempty.
The common intersection of a family fixed by all left shifts in an ordered set is an initial interval there.
Apply this observation to the other single-index section to obtain the final assertion in \textup{(iii)}.
If $|C\cup D|\leq k-1$, local shiftedness puts in both sections the set consisting of $C\cup D$
together with the first $k-1-|C\cup D|$ points outside it, again a contradiction.

Finally, suppose $H\in\mathcal H$ omits $c\in C\setminus D$ and $d\in D\setminus C$.
The corresponding member with both critical coordinates avoids $c,d$, so their left shift in the original order is fixed.
If $c<d$, apply it to a member of the second section omitting $c$; it removes the compulsory point $d$.
If $d<c$, use the first section instead. Both conclusions are impossible.
This proves \textup{(iv)}.
\end{proof}

The next two lemmas exclude the small-core configurations that arise
in the proof of Theorem~\ref{mth}.

\begin{lem}\label{mth-star-estimates}
Let $n>2k$, $k\ge3$, and $a\ge2$. Then the following hold.
\begin{enumerate}
\item[\textup{(i)}]
If $P\cap L=\emptyset$, $|L|=k$, and $p=|P|\in\{2,3\}$, then
\[
\mathcal F=\{F\in\binom{[n]}k:P\subseteq F\}\cup\{L\}
\]
satisfies $a|\mathcal F|+|\overline N_k(\mathcal F)|<(a+1)h(n,k,1)$.
\item[\textup{(ii)}]
Let $k\ge4$, and let $C,D\in\binom{[3,n]}{k-1}$ be distinct, with $C\cap D\ne\emptyset$. If $x\in[3,n]\setminus(C\cap D)$, then
\[
\mathcal F=\{\{1\}\cup C,\{2\}\cup D\} \cup\{\{1,2\}\cup H:H\in\binom{[3,n]}{k-2},\ x\in H\}
\]
satisfies $a|\mathcal F|+|\overline N_k(\mathcal F)|<(a+1)h(n,k,1)$.
\end{enumerate}
\end{lem}

\begin{proof}
For \textup{(i)}, a $k$-set meets every $k$-set containing $P$ if and only if it meets $P$, since $n>2k$. Thus
\[
|\mathcal F|=\binom{n-p}{k-p}+1, \qquad |\overline N_k(\mathcal F)| =\sum_{i=1}^{p}\left(\binom{n-i}{k-1}-\binom{n-k-i}{k-1}\right).
\]
Since $p\ge2$,
\[
|\mathcal F|\le\binom{n-2}{k-2}+1<h(n,k,1),
\]
where the strictness follows from
$\binom{n-2}{k-1}>\binom{n-k-1}{k-1}$.
The difference $(a+1)h(n,k,1)-a|\mathcal F|-|\overline N_k(\mathcal F)|$
therefore increases with $a$, so it suffices to take $a=2$.
For $p=2$,
\[
3h(n,k,1)-2|\mathcal F|-|\overline N_k(\mathcal F)|
=\sum_{i=3}^{k+1}\binom{n-i}{k-2}
-\binom{n-k-2}{k-2}+1>0.
\]
If $p=3$, the same difference equals
\[
3\binom{n-3}{k-2}-2\binom{n-k-2}{k-2} -\binom{n-k-3}{k-2}+1>0.
\]
This proves \textup{(i)}.

For \textup{(ii)}, again
$|\mathcal F|=\binom{n-3}{k-3}+2<h(n,k,1)$.
The companion's double section has size $\binom{n-2}{k-2}$,
and each single-index section has size
$\binom{n-2}{k-1}-\binom{n-k-1}{k-1}$.
A member avoiding $1,2$ must contain $x$ and meet both $C,D$.
At least one of these sets avoids $x$, so this section has size
at most $\binom{n-3}{k-1}-\binom{n-k-2}{k-1}$.
It follows that
\[
3h(n,k,1)-2|\mathcal F|-|\overline N_k(\mathcal F)| \ge3\binom{n-3}{k-2}-\binom{n-k-2}{k-2}-1>0.
\]
Since $|\mathcal F|<h(n,k,1)$, the result follows for all $a\ge2$.
\end{proof}

\begin{lem}\label{mth-generated-exceptions}
Let $n,k,a$ be integers with $n>2k$, $k\geq3$, $a\geq2$, and
\[
a\binom{n-2}{k-2}<\binom{n-2}k,\qquad
M=\max\{(a+1)h(n,k,1),g(n,k)+3a\}.
\]
Suppose $\mathcal P$ is a nonempty family of $k$-sets with $\bigcap\mathcal P=P$.
\begin{enumerate}
\item[\textup{(i)}] If $|P|=2$, $|L|=k$, and $P\cap L=\emptyset$, then
$\mathcal F=\mathcal P\cup\{L\}$ satisfies
$a|\mathcal F|+|\overline N_k(\mathcal F)|<M$.
If $k\geq4$, the stronger inequality with an additional $a$ on the left also holds.
\item[\textup{(ii)}] If $|P|=3$, $k\geq4$, $|L|=k$, $P\cap L=\emptyset$, and $|\mathcal P|\geq3$,
then $\mathcal F=\mathcal P\cup\{L\}$ satisfies the same strict bound.
\item[\textup{(iii)}] If $|P|=3$, $k\geq4$, $|Q|=k-1$, $P\cap Q=\emptyset$, and $x,y$ are distinct points of $P$,
then $\mathcal F=\mathcal P\cup\{Q\cup\{x\},Q\cup\{y\}\}$ satisfies the same strict bound.
\end{enumerate}
\end{lem}

\begin{proof}
Put $v=n-k$. Removing the common core $P$ gives the nonempty nontrivial family
\[
\mathcal R=\{F\setminus P:F\in\mathcal P\}
 \subseteq\binom{[n]\setminus P}{k-|P|}.
\]
For $p=|P|\in\{2,3\}$, the coefficient condition needed in Theorem~\ref{m2n} holds, since
\[
\frac{\binom{n-p}k}{\binom{n-p}{k-p}}
=\frac{\binom{n-2}k}{\binom{n-2}{k-2}}
  \prod_{i=2}^{p-1}\frac{n-k-i}{k-i}>a.
\]

For \textup{(i)} and \textup{(ii)}, companion members meeting $P$ must also meet $L$;
there are exactly
\[
Q_p=\binom nk-\binom{n-p}k-\binom vk+\binom{v-p}k
\]
such members. Those avoiding $P$ are contained in the full companion of $\mathcal R$.
Theorem~\ref{m2n} therefore gives
\begin{equation}\label{mth-one-exception-endpoints}
a|\mathcal F|+|\overline N_k(\mathcal F)|
\leq\max\{U_p,\,G_p+3a\},
\end{equation}
where
\[
U_p=a\left(\binom{n-p}{k-p}+1\right)+Q_p,\qquad
G_p=\binom nk-3\binom vk+\binom{v-p}k+\binom{v-k+p}k.
\]
The value $U_p$ corresponds to the full principal family generated
by $P$, together with $L$.
Lemma~\ref{mth-star-estimates}\textup{(i)} gives $U_p<(a+1)h(n,k,1)$.
Discrete convexity gives $G_p\leq g(n,k)$ for $1\leq p\leq k-1$.
For $p=2$, the comparison is strict if $k\geq4$; for $k=3$,
\[
(a+1)h(n,3,1)-g(n,3)-3a=(a-2)(3n-11)+2n-5>0.
\]
This proves the first assertion in \textup{(i)}.
For $p=3$, the comparison is strict if $k\geq5$.
When $k=4$, equality in the potentially maximal second endpoint of
\eqref{mth-one-exception-endpoints} would, by Theorem~\ref{m2n}, require $|\mathcal R|=2$;
the full-level equality alternative gives $U_p<M$ instead.
The hypothesis $|\mathcal P|\geq3$ excludes this possibility, proving \textup{(ii)}.

For the stronger part of \textup{(i)}, observe first that
\begin{equation}\label{mth-strong-coefficient}
a<\frac{v(v-1)}{k(k-1)}<\binom{v-3}{k-2}\qquad(k\geq4).
\end{equation}
Indeed, $a\geq2$ forces $v\geq k+2$.
The second inequality holds at $v=k+2$, and the ratio
$\binom{v-3}{k-2}/(v(v-1))$ strictly increases: its successive ratio is
$(v-2)(v-1)/((v-k)(v+1))>1$.
Consequently
\[
g(n,k)-G_2
=\binom{v-2}{k-1}-\binom{v-k+1}{k-1}
\geq\binom{v-3}{k-2}>a.
\]
Also $U_2+a<(a+1)h(n,k,1)$.
For this last assertion it is enough to take $a=2$, since
$\binom{n-2}{k-2}+2<h(n,k,1)$; the difference at $a=2$ is
\[
\sum_{i=3}^{k+1}\binom{n-i}{k-2}-\binom{n-k-2}{k-2}-1>0.
\]
Both endpoints in \eqref{mth-one-exception-endpoints}, even after adding $a$, are below $M$.

For \textup{(iii)}, let $p=3$.
The number of companion members meeting $P$ is exactly
\[
Z=\binom nk-2\binom vk+\binom{v-1}k
       -\binom{n-3}k+\binom{v-2}k.
\]
To obtain this expression, count the full companion of the two exceptional members,
and subtract the $k$-sets avoiding $P$ and meeting $Q$.
The remaining companion members are again contained in the full companion of $\mathcal R$.
Theorem~\ref{m2n} bounds the objective by $\max\{V,W\}$, where
\[
V=a\left(\binom{n-3}{k-3}+2\right)+Z,
\qquad
W=Z+\binom{n-3}k-2\binom vk+\binom{v-k+3}k+4a.
\]
For the first value, $V<(a+1)h(n,k,1)$:
the difference at $a=2$ is
\[
3\binom{n-3}{k-2}-\binom{v-2}{k-2}-1>0,
\]
and $\binom{n-3}{k-3}+2<h(n,k,1)$.
Finally,
\[
g(n,k)+3a-W
=\binom{v-1}{k-1}+\binom{v-2}{k-1}
-\binom{v-k+2}{k-1}-\binom{v-k+1}{k-1}-a
\geq\binom{v-2}{k-2}-a>0.
\]
The last inequality follows from \eqref{mth-strong-coefficient}.
Thus $W<M$ as well, proving \textup{(iii)}.
\end{proof}

\begin{proof}[Proof of Theorem~\ref{mth}]
The case $k=1$ is vacuous.
For $k=2$, if all families are intersecting, nontriviality makes
each a triangle, and cross-intersection makes the triangles coincide.
Otherwise, one family contains disjoint edges.
Every other family then lies in the bipartite graph of edges
meeting both of them. It cannot be nontrivial and intersecting,
so it too contains disjoint edges.
Applying Lemma~\ref{abc} to two of the other families bounds the
first family's size by $g(n,2)-1=2$.
Thus every family has size $2$, and the total is less than $3m$.
This proves both assertions for $k=2$.

Assume that $k\ge3$.
Put $a=m-1$ and $M=\max\{(a+1)h(n,k,1),g(n,k)+3a\}$.
Choose nonempty, nontrivial cross-intersecting families
$\mathcal A_1,\ldots,\mathcal A_m\subseteq\binom{[n]}k$
maximizing $\Sigma=\sum_{i=1}^m|\mathcal A_i|$, ordered by
$|\mathcal A_1|\geq\cdots\geq|\mathcal A_m|$.
The two constructions give $\Sigma\geq M$.
Hence $|\mathcal A_1|\geq h(n,k,1)$, and Theorem~\ref{mh} gives
$|\mathcal A_2|\leq h(n,k,1)\leq|\mathcal A_1|$.

If $|\mathcal A_2|=2$, each $\mathcal A_i$ with $i\geq2$
consists of two disjoint $k$-sets.
Lemma~\ref{abc}, applied to $\mathcal A_2,\mathcal A_3$, yields
$|\mathcal A_1|<g(n,k)$ and therefore
$\Sigma<g(n,k)+2a<M$, a contradiction.
Set $\mathcal A=\mathcal A_2$ and $\mathcal B=\mathcal A_1$.
Then
\begin{equation}
\label{mth-initial-comparison} 3\leq|\mathcal A|\leq h(n,k,1)\leq|\mathcal B|, \qquad W:=a|\mathcal A|+|\mathcal B|\geq\Sigma\geq M.
\end{equation}

We will use the following consequence of extremality. If $\mathcal F$ is nontrivial and intersecting, and $\mathcal G$ is nontrivial and
cross-intersecting with $\mathcal F$, then $\mathcal G,\mathcal F,\ldots,\mathcal F$ is an admissible $m$-tuple. Therefore
\begin{equation}
\label{mth-tuple-comparison} a|\mathcal F|+|\mathcal G|\leq\Sigma.
\end{equation}

A simultaneous shift preserving nontriviality of $\mathcal A$
also preserves that of $\mathcal B$.
Indeed, if the shifted $\mathcal B$ had a common point $x$,
the shifted $\mathcal A$ would have a member avoiding $x$.
Cross-intersection would then give
$|\mathcal B|\leq\binom{n-1}{k-1}-\binom{n-k-1}{k-1}
=h(n,k,1)-1$, a contradiction.
Apply such shifts until $\mathcal A$ is shifted or critically shifted.

If $\mathcal A$ is shifted, Theorem~\ref{sab} bounds $W$ by
the larger of $(a+1)h(n,k,1)$ and the clique value.
Lemma~\ref{mklc} makes the latter strictly smaller than $M$.
Since $W\geq M$, we obtain
$\Sigma\leq W\leq(a+1)h(n,k,1)$ and
$M=(a+1)h(n,k,1)$.
Equality makes the resulting pair a common maximum nontrivial
intersecting family. We recover the original families at the end.

Suppose now that $\mathcal A$ is critically shifted.
Enlarge $\mathcal B$ to $\overline N_k(\mathcal A)$.
This preserves nontriviality and can only increase $W$.

Apply Lemma~\ref{critical-local-generators} to $\mathcal A$ and use its ordered remaining coordinates as $[3,n]$.
Set $C=\bigcap\mathcal A(1)$, $D=\bigcap\mathcal A(2)$, $\mathcal H=\mathcal A(12)$,
$c=|C|$, $d=|D|$, $T=\binom{n-2}{k-2}$, and $\mathcal B_0=\mathcal B(\emptyset)$.
We have
\[
\mathcal A(\emptyset)=\emptyset,\quad
\mathcal A(1)\cap\mathcal A(2)=\emptyset,\quad
\mathcal B(12)=\binom{[3,n]}{k-2},\quad
C,D\neq\emptyset,\quad |C\cup D|\geq k.
\]
In particular, $C\setminus D$ and $D\setminus C$ are nonempty.
The deletion property is
\begin{equation}\label{mth-deletion}
E\setminus\{j\}\in\mathcal H
\quad(E\in\mathcal A(1),\ j\in E\setminus C),
\end{equation}
and likewise for the second section.
Also,
\begin{equation}\label{mth-local-core-cover}
H\supseteq C\setminus D\quad\hbox{or}\quad H\supseteq D\setminus C
\qquad(H\in\mathcal H).
\end{equation}

We have $0<|\mathcal H|<T$.
Indeed, if $\mathcal H=\emptyset$, deletion makes both single-index
sections singletons, contradicting $|\mathcal A|\geq3$.
For the upper bound, choose one point from each of
$C\setminus D$ and $D\setminus C$.
A $(k-2)$-set avoiding both exists and is excluded by
\eqref{mth-local-core-cover}.

Extremality gives the further comparison
\begin{equation}\label{mth-core-optimality}
a|\mathcal H|+|\mathcal B_0|\geq aT.
\end{equation}
Replace $\mathcal H$ by its full level and delete $\mathcal B_0$.
The new pair remains cross-intersecting. Shift simultaneously within $[3,n]$ until both single-index sections
of the new $\mathcal A$ are shifted. Shifting preserves triviality of a nonempty family:
if a common point moves, every member acquires the target point.
Thus both sections have $3$ as a common point, and the new $\mathcal A$ is intersecting.
Its full double section and nonempty single sections ensure nontriviality.
The new $\mathcal B$ is also nontrivial: its double section is full, and its two single sections are nonempty,
since before shifting they contain all $(k-1)$-sets meeting $D$ and $C$, respectively.
Applying \eqref{mth-tuple-comparison} to this new tuple proves
\eqref{mth-core-optimality}. Since $|\mathcal H|<T$, it also gives
$\mathcal B_0\ne\emptyset$.

We now use this comparison to bound $a$.
The disjointness graph between levels $k-2$ and $k$ of $[3,n]$
is biregular and connected.
Indeed, any two $(k-2)$-sets differing in one point have a common disjoint $k$-set, since $n>2k$;
such one-point changes connect that level, and every vertex of the other level has a neighbor.
Counting edges out of $\mathcal H$ therefore gives
$|N_k(\mathcal H)|/\binom{n-2}k\geq|\mathcal H|/T$.
Equality would separate $\mathcal H\cup N_k(\mathcal H)$ from its complement,
contradicting connectedness because $0<|\mathcal H|<T$.
Since $\mathcal B_0$ avoids this neighborhood, we obtain
\[
\frac{|\mathcal H|}{T}+\frac{|\mathcal B_0|}{\binom{n-2}k}<1.
\]
Together with \eqref{mth-core-optimality}, this forces
\begin{equation}\label{mth-coefficient}
a\binom{n-2}{k-2}<\binom{n-2}k.
\end{equation}

Put $K=\bigcap\mathcal H$. Whenever $c\leq k-2$, choose $E\in\mathcal A(1)$.
Then $E\setminus C\neq\emptyset$, and \eqref{mth-deletion} gives
\begin{equation}\label{mth-core-contained}
K\subseteq\bigcap_{j\in E\setminus C}(E\setminus\{j\})=C.
\end{equation}
Likewise, $K\subseteq D$ whenever $d\leq k-2$.

\medskip
\noindent\textit{Claim.} Neither $c$ nor $d$ lies between $3$ and $k-2$.

Suppose that $3\leq c\leq k-2$.
Under this hypothesis, we first establish
\begin{equation}\label{mth-forced-star}
\{H\in\tbinom{[3,n]}{k-2}:C\subseteq H\}\subseteq\mathcal H,
\qquad \bigcap\mathcal H\subseteq C.
\end{equation}
The second inclusion is \eqref{mth-core-contained}.
To prove the first, suppose that some $(k-2)$-set containing $C$ is missing from $\mathcal H$.
We will construct an admissible generator contraction with a
strictly larger weighted sum.
Lemma~\ref{critical-local-generators} supplies
$E_0\in\mathcal A(1)$ meeting $D$.
If $C\cap D\neq\emptyset$, every member does so.
Otherwise, $D\setminus C$ is a nonempty initial interval of the residual ground set,
and the shifted residual family contains its first $(k-c-1)$-set, which meets that interval.
Choose a canonical generator $G\subseteq E_0$ of $\mathcal A(1)$
on $[3,n]$. Delete the other generators and enlarge $G$ to $E_0$,
contracting the section to $\{E_0\}$.
Let $\mathcal A'$ be the resulting family.
Replace $\mathcal B(2)$ by all $(k-1)$-sets in $[3,n]$ meeting
$E_0$ to obtain $\mathcal B'$, and keep every other section fixed.

The deletion property \eqref{mth-deletion} gives
\[
\{E_0\setminus\{j\}:j\in E_0\setminus C\}\subseteq\mathcal H,
\]
so the shadow constraint is preserved with $\mathcal H$ unchanged.
New members of the $2$-section of $\mathcal B'$ meet $\{1\}\cup E_0$ and share $2$
with the other members of $\mathcal A'$. All other cross-intersections are preserved by
$\mathcal A'\subseteq\mathcal A$.
Moreover, $\mathcal A'$ is intersecting because $E_0\cap D\neq\emptyset$.
The inclusion $K\subseteq C\subseteq E_0$ and nontriviality
of the original family give
\[
\bigcap\mathcal A'=E_0\cap D\cap K=D\cap K=C\cap D\cap K=\emptyset.
\]
Hence $\mathcal A'$ remains nontrivial.
Finally, $\mathcal B'\supseteq\mathcal B$ is nontrivial.
Thus $\mathcal B',\mathcal A',\ldots,\mathcal A'$ is an admissible $m$-tuple.

We now count the new members allowed in $\mathcal B(2)$.
Put $N=n-c-2$, $p=k-c-1$, and $q=k-1$.
For this count, identify $[3,n]\setminus C$ with $[N]$ so that $E_0\setminus C=[p]$,
and write $\mathcal R=\{F\setminus C:F\in\mathcal A(1)\}$ in these coordinates.
For $1\leq i\leq p$, set $V_i=[i+1,N]$ and
\[
\mathcal R_i=
\{R\setminus[i-1]:R\in\mathcal R\setminus\{[p]\},\ \min([p]\setminus R)=i\}.
\]
Each $(q-1)$-set $Q\subseteq V_i$ disjoint from some member of $\mathcal R_i$
produces a newly available member $\{i\}\cup Q$ of $\mathcal B'(2)$.
Indeed, it meets $E_0$ at $i$ and is disjoint from the old member $C\cup[i-1]\cup R$
corresponding to $R\in\mathcal R_i$: $Q$ avoids $R\cup[i]$, and $i\notin R$.
These newly available families are disjoint for different $i$, since their least coordinates differ.
Counting disjoint pairs on $V_i$ gives at least
\[
\frac{\binom{N-p-1}{q-1}}{\binom{N-i-q+1}{p-i+1}}|\mathcal R_i|
=\frac{\binom{N-i}{q-1}}{\binom{N-i}{p-i+1}}|\mathcal R_i|
\]
new members for index $i$: these are the degree ratios of the two levels of the disjointness graph.
The last ratio increases with $i$, since its successive ratio is
$(N-i-q+1)/(p-i+1)>1$.
As $\sum_i|\mathcal R_i|=|\mathcal A(1)|-1$, we obtain
\begin{equation}\label{mth-generator-gain}
|\mathcal B'(2)|-|\mathcal B(2)|
\geq\Lambda_c\bigl(|\mathcal A(1)|-1\bigr),\qquad
\Lambda_c=\frac{\binom{n-c-3}{k-2}}{\binom{n-c-3}{k-c-1}}.
\end{equation}

Put $v=n-k>k$. The quantities $\Lambda_c$ increase with $c$, since
$\Lambda_{c+1}/\Lambda_c=(v-c-1)/(k-c-1)>1$.
Consequently,
\[
\Lambda_c\geq\Lambda_3
=\frac{(v-2)(v-3)}{(k-2)(k-3)}
>\frac{v(v-1)}{k(k-1)}>a.
\]
The middle inequality follows by comparing $(v-2)/(k-2)$ with $v/k$ and
$(v-3)/(k-3)$ with $(v-1)/(k-1)$; both comparisons are strict because $v>k$.
The last inequality is exactly \eqref{mth-coefficient}.
Since $c<k-1$ implies $|\mathcal A(1)|>1$, equation~\eqref{mth-generator-gain} gives
\begin{equation}\label{mth-generator-improvement}
a|\mathcal A'|+|\mathcal B'|-W
\geq(\Lambda_c-a)\bigl(|\mathcal A(1)|-1\bigr)>0.
\end{equation}
This contradicts \eqref{mth-tuple-comparison}, since $W\geq\Sigma$.
Hence no $(k-2)$-set containing $C$ can be missing from $\mathcal H$, proving \eqref{mth-forced-star}.

The same contraction now completes the claim.
Its admissibility follows from deletion and \eqref{mth-core-contained};
the strict gain uses only $3\leq c\leq k-2$ and
\eqref{mth-coefficient}.
These conditions still hold when the full $C$-star is present,
so \eqref{mth-generator-improvement} excludes that case as well.
Thus $3\leq c\leq k-2$ is impossible.
The argument for $d$ is symmetric, and hence
\[
c,d\in\{1,2,k-1\}.
\]

It remains to exclude the small core sizes using
Lemma~\ref{mth-generated-exceptions}.
If $c=1$, then $d=k-1$ and $C\cap D=\emptyset$.
Equation~\eqref{mth-local-core-cover} forces every member of $\mathcal H$ to contain $C$.
All members except $\{2\}\cup D$ therefore contain $P=\{1\}\cup C$, and their common intersection is exactly $P$.
The exceptional member is disjoint from $P$.
Lemma~\ref{mth-generated-exceptions}\textup{(i)} contradicts $W\geq M$.
Treat $d=1$ symmetrically.

Suppose next that $c=2<k-1$.
Unless $k=4$ and $d=2$, we have $d=k-1$.
If $C\cap D=\emptyset$, equation~\eqref{mth-local-core-cover} makes every member of $\mathcal H$ contain $C$.
The part containing $P=\{1\}\cup C$ has core exactly $P$ and at least three members:
the first section has at least two, and deletion supplies a member in the double section.
The only other member is $\{2\}\cup D$, disjoint from $P$.
Lemma~\ref{mth-generated-exceptions}\textup{(ii)} gives a contradiction.

If $|C\cap D|=1$, write $C=\{x,y\}$ with $y\in D$ and $x\notin D$.
Every double-section member contains $x$ or equals $D\setminus C$.
If an $x$-containing double-section member avoids $y$, then, after removing the possible member
$\{1,2\}\cup(D\setminus C)$, the family consists of a part with core $\{1,x\}$ and the disjoint exceptional member
$\{2\}\cup D$. The stronger assertion in
Lemma~\ref{mth-generated-exceptions}\textup{(i)}
absorbs the contribution $a$ of the removed member and still gives
$W<M$.
Otherwise all $x$-containing double-section members contain $C$.
Nontriviality forces $D\setminus C\in\mathcal H$.
The family now consists of a part with core $P=\{1,x,y\}$ and the two exceptional members
\[
Q\cup\{y\},\quad Q\cup\{1\},\qquad Q=\{2\}\cup(D\setminus C),
\]
where $P\cap Q=\emptyset$ and $|Q|=k-1$.
Lemma~\ref{mth-generated-exceptions}\textup{(iii)} again gives $W<M$.
These are all possibilities, since $C\subseteq D$ would contradict $|C\cup D|\geq k$.

It remains to treat $k=4$ and $c=d=2$.
Here $C\cap D=\emptyset$ and deletion and \eqref{mth-local-core-cover} give $\mathcal H=\{C,D\}$.
Put $R=\{1\}\cup C$ and $S=\{2\}\cup D$.
The family has the form
\[
\{R\cup\{x\}:x\in X\}\ \dot\cup\ \{S\cup\{y\}:y\in Y\},
\]
where $X\subseteq[n]\setminus R$ and $Y\subseteq[n]\setminus S$.
Local shiftedness implies $S\subseteq X$ and $R\subseteq Y$:
the residual single-index families are initial intervals of size at least two, and their first two points are
$D$ and $C$, respectively. Thus $3\leq|X|,|Y|\leq n-3$.
A companion member meets $R$ or contains $X$, and meets $S$ or contains $Y$.
The generator rule gives
\[
W= a(|X|+|Y|)+Z+f(|X|)+f(|Y|),\qquad
Z=\binom n4-2\binom{n-3}4+\binom{n-6}4,
\]
where $f(3)=n-6$, $f(4)=1$, and $f(t)=0$ for $t\geq5$.
For every $3\leq t\leq n-3$, $at+f(t)\leq a(n-3)$, since $a\geq2$ and $n\geq9$.
Consequently $W\leq2a(n-3)+Z<(a+1)h(n,4,1)$.
For the last inequality it is enough to take $a=2$; the gap is
$(3n-14)(n-3)/2>0$, and $2(n-3)<h(n,4,1)$.
All small-core cases are therefore excluded, and we have
\begin{equation}\label{mth-singleton-sections}
\mathcal A(1)=\{C\},\qquad \mathcal A(2)=\{D\},\qquad |C|=|D|=k-1.
\end{equation}
Since $\mathcal H\ne\emptyset$, \eqref{mth-local-core-cover}
also gives $C\cap D\ne\emptyset$.
Thus $\mathcal A$ is intersecting.

For $k=3$, the remaining possibilities can be counted directly.
In this case $|C\cap D|=1$, and $\mathcal H$ is a nonempty subfamily of the two singletons
$C\setminus D,D\setminus C$.
If $|\mathcal H|=1$, the section counts give $|\mathcal B|\leq g(n,3)=7n-25$,
so $W\leq3a+g(n,3)<(a+1)h(n,3,1)$.
If $|\mathcal H|=2$, the double, two single, and empty sections of $\mathcal B$ have respective sizes
$n-2,2n-7,2n-7,n-4$. Hence
\[
W=4a+6n-20<(a+1)(3n-8),
\]
with gap $(a-1)(3n-12)>0$.
Both are contradictions. Henceforth $k\geq4$.

All remaining neighborhoods are taken on $[3,n]$.
Put $v=n-k$.
We first show that the double section must be a singleton.
Suppose that $|\mathcal H|\ge2$, and put
$K=\bigcap\mathcal H$, $t=|K|+2$, and $q=k-t$.
We need the estimate
\begin{equation}
\label{mth-neighborhood-lower} |N_k(\mathcal H)|-a|\mathcal H| \geq2\binom vk-\binom{v-q}k-2a.
\end{equation}
The family $\mathcal H'=\{H\setminus K:H\in\mathcal H\}$ is nonempty and nontrivial in $\binom{[3,n]\setminus K}q$. Also,
\eqref{mth-coefficient} implies
\[
\frac{\binom{n-t}k}{\binom{n-t}q} = \frac{\binom{n-2}k}{T}\prod_{j=0}^{t-3}\frac{n-k-2-j}{k-2-j}>a.
\]
Theorem~\ref{m2n}, applied to $\mathcal H'$ and its full cross-intersecting $k$-uniform family on $[3,n]\setminus K$, gives
\begin{equation}
\label{mth-core-weighted}
a|\mathcal H|+\binom{n-t}k-|N_k(\mathcal H)|
\le\max\left\{\binom{n-t}k-2\binom vk+\binom{v-q}k+2a,\,
a\binom{n-t}q\right\}.
\end{equation}
We show that the first expression in this maximum is at least
the second. For $t\ge4$, apply Lemma~\ref{m2n-binomial}\textup{(ii)}
with parameters $(n-2,k-2,k,t-2)$ and coefficient $a$.

For $t=2$, the second expression is $aT$.
If it were larger, the equality statement in Theorem~\ref{m2n}
and $|\mathcal H|<T$ would make the bound $aT$ in
\eqref{mth-core-weighted} strict.
Since $\mathcal B_0$ lies in the full companion of $\mathcal H$,
this would contradict \eqref{mth-core-optimality}.

For $t=3$, write $K=\{x\}$ and suppose that the second expression
is larger. Separate $\mathcal B_0$ according to whether its
members contain $x$.
Theorem~\ref{m2n}, applied to the members avoiding $x$ and to
$\mathcal H'$, bounds their weighted sum by $a\binom{n-t}q$.
We may therefore replace $\mathcal H$ by the full $x$-star
and $\mathcal B_0$ by all $k$-sets in $[3,n]$ containing $x$
and meeting both $C,D$, without decreasing the objective.
The other sections remain fixed.
The new pair is cross-intersecting and nontrivial.
Since $x\notin C\cap D$, Lemma~\ref{mth-star-estimates}\textup{(ii)}
bounds its weighted sum strictly below $(a+1)h(n,k,1)\le M$,
a contradiction.
This proves \eqref{mth-neighborhood-lower}.

Put $\alpha=|C\setminus D|+1$, $J=|N_k(\mathcal H)\cap N_k(\{C,D\})|$. Since $\mathcal B$ is saturated,
\[
|\mathcal B| = \binom nk-2\binom vk+\binom{v-\alpha}k -|N_k(\mathcal H)|+J.
\]
Define
\[
\Delta=|N_k(\mathcal H)|-a|\mathcal H|-J-\binom{v-\alpha}k+\binom{v-k}k
\]
and $\Gamma=\binom{v-1}{k-1}-\binom{v-k}{k-1}$. Since $|\mathcal A|=|\mathcal H|+2$, we have
\begin{equation}
\label{mth-objective-identity} W=g(n,k)+3a+\Gamma-a-\Delta.
\end{equation}

Let $K=\bigcap\mathcal H$, and put
\[
x=|K\cap(C\setminus D)|+1,\quad y=|K\cap(D\setminus C)|+1,\quad z=|K\setminus(C\cup D)|,\quad t=|K|+2.
\]
Nontriviality gives $K\cap C\cap D=\emptyset$, so $x+y+z=t$ and $1\leq x,y\leq\alpha$. Every set counted by $J$ avoids $K$ and at least one
of $C,D$. Consequently, $J\leq\binom{v-x-z}k+\binom{v-y-z}k-\binom{v-\alpha-z}k$. Discrete convexity gives
\begin{align*}
\binom{v-\alpha}k-\binom{v-\alpha-z}k &\leq\binom{v-x}k-\binom{v-x-z}k, \\
\binom{v-x}k+\binom{v-y-z}k &\leq\binom{v-1}k+\binom{v-t+1}k.
\end{align*}
Thus
\begin{equation}
\label{mth-J-bound} J+\binom{v-\alpha}k \leq\binom{v-1}k+\binom{v-t+1}k.
\end{equation}

If $|\mathcal H|\geq2$, combine
\eqref{mth-neighborhood-lower} and \eqref{mth-J-bound}.
Since $q+t=k$, this gives
\[
\Delta \geq\Gamma+ \sum_{i=1}^{q}\sum_{j=1}^{t-1}\binom{v-i-j}{k-2}-2a \geq\Gamma+\binom{v-2}{k-2}-2a>\Gamma-a.
\]
The last inequality follows from Lemma~\ref{m2n-binomial}(i) with $(n,k,\ell)$ replaced by $(n-2,k-2,k)$.
Equation~\eqref{mth-objective-identity} gives $W<g(n,k)+3a$, a contradiction. Hence $\mathcal H=\{K\}$, where $|K|=k-2$. Now $t=k$ and
$|N_k(\mathcal H)|=\binom vk$. Equation~\eqref{mth-J-bound} gives $\Delta\geq\Gamma-a$, and therefore
\begin{equation}
\label{mth-critical-bound} W\leq g(n,k)+3a.
\end{equation}

For $\mathcal H=\{K\}$, the initial bound for $J$ is exact.
Equality in \eqref{mth-J-bound} holds precisely when $z=0$
and $x=1$ or $y=1$.
Indeed, convexity is strict when $x>1$ and $y+z>1$;
if $x=1$ and $z>0$, the preceding binomial-difference
inequality is strict.
Here $\alpha\geq2$ and $\binom{v-3}{k-2}>0$.
Thus equality requires
$K\subseteq C\setminus D$ or $K\subseteq D\setminus C$.
Since $C\cap D\ne\emptyset$ and $|K|=k-2$, we obtain
$|C\cap D|=1$ and $K=C\setminus D$ or $K=D\setminus C$.
Up to isomorphism, this gives
\begin{equation}
\label{mth-extremal-triple} \mathcal A = \{\{1\}\cup[3,k+1],\{2\}\cup[k+1,2k-1],[k]\}, \qquad \mathcal B=\overline N_k(\mathcal A).
\end{equation}

It remains to recover the original equality families.
First consider equality in the shifted case.
All original families have size $h(n,k,1)$.
Apply the same shifts to all $m$ families.
The size argument above preserves their nontriviality.
At the end, every family is cross-intersecting with the same
maximum nontrivial intersecting family and has the same size,
so all families coincide.

We recover this property through each inverse shift.
Suppose $\mathcal F_1,\ldots,\mathcal F_m$ have the same
intersecting image under $s_{ij}$.
If $\mathcal F_1$ had disjoint members, they would be
$\{i\}\cup X$ and $\{j\}\cup Y$, with
$X,Y\subseteq[n]\setminus\{i,j\}$ and $X\cap Y=\emptyset$.
The common image contains $\{i\}\cup X$ and $\{i\}\cup Y$.
Cross-intersection forces each $\mathcal F_r$, $r\geq2$,
to contain $\{j\}\cup X$ and $\{i\}\cup Y$.
Taking these members from $\mathcal F_2$ and $\mathcal F_3$
contradicts cross-intersection.
Thus every original family is intersecting.
Members of one can then be added to the others, so extremality
and equal sizes force the families to coincide.
Repeating this argument recovers case \textup{(i)}.

Now consider equality in \eqref{mth-critical-bound}.
All preceding size comparisons are equalities, so the original
families satisfy $|\mathcal A_1|=g(n,k)$ and
$|\mathcal A_i|=3$ for $i\geq2$.
Two smaller families cannot both contain disjoint pairs:
Lemma~\ref{abc} would give $|\mathcal A_1|<g(n,k)$.
Hence at least one smaller family, say $\mathcal F$, is intersecting.
Its members can be added to every other smaller family without
destroying admissibility.
Extremality and equal sizes give
$\mathcal A_2=\cdots=\mathcal A_m=\mathcal F$ and
$\mathcal A_1=\overline N_k(\mathcal F)$.

The original $\mathcal F$ is a nontrivial intersecting triple, and the shifts carry it to \eqref{mth-extremal-triple}.
Track the three members through each shift $s_{ij}$. A pairwise intersection can decrease only if two members contain $j$ and avoid $i$,
one of them moves and the other is retained. Write them as $\{j\}\cup X$ and $\{j\}\cup Y$, with the former retained.
Retention forces $\{i\}\cup X$ to be the third member. The common intersection of the triple is then $X\cap Y$,
so nontriviality gives $X\cap Y=\emptyset$. But this makes $\{j\}\cup Y$ disjoint from $\{i\}\cup X$, a contradiction.
Thus no pairwise intersection size decreases along the shifts.
The final intersection sizes are $(k-1,1,1)$, and the original triple is intersecting; its intersection sizes are therefore
$(s,1,1)$ for some $1\leq s\leq k-1$. Hence
\[
|\overline N_k(\mathcal F)|= \binom nk-3\binom{n-k}k+\binom{n-2k+s}k +2\binom{n-2k+1}k-\binom{n-3k+s+2}k.
\]
This expression equals $g(n,k)$ when $s=k-1$.
Its increment from $s$ to $s+1$ is
\[
\binom{n-2k+s}{k-1}-\binom{n-3k+s+2}{k-1}\geq0.
\]
The increment from $s=k-2$ to $s=k-1$ is strictly positive for $n>2k$ and $k\ge4$. Since $|\overline N_k(\mathcal F)|=g(n,k)$, we obtain
$s=k-1$. The original triple is therefore isomorphic to
\eqref{mth-extremal-triple}, giving case (ii).
The two constructions give the converse and complete the proof.
\end{proof}

\end{document}